\documentclass{aims}
\usepackage{amsmath}
\usepackage{paralist}
\usepackage{graphics} 
\usepackage{epsfig} 
\usepackage{graphicx}
\usepackage{epstopdf}
\usepackage[colorlinks=true]{hyperref}
\hypersetup{urlcolor=blue, citecolor=red}
\usepackage{subcaption}

\def\currentvolume{X}
 \def\currentissue{X}
  \def\currentyear{200X}
   \def\currentmonth{XX}
    \def\ppages{X--XX}

\newtheorem{theorem}{Theorem}[section]

\newtheorem{lemma}[theorem]{Lemma}
\newtheorem{proposition}{Proposition}

\theoremstyle{definition}

\title[COMPETING SPECIES SHARING A PARASITE]
      {A DISCRETE MODEL OF COMPETING SPECIES \newline SHARING A PARASITE}

\author[RAFAEL BRAVO DE LA PARRA and LUIS SANZ]{}



 \email{rafael.bravo@uah.es}
 \email{luis.sanz@upm.es}

\thanks{Authors are supported by Ministerio de Econom\'{\i}a y Competitividad (Spain),
project MTM2014-56022-C2-1-P.}

\begin{document}

{\centering%
\noindent\begin{tabular}{|l|}
\hline
\textit{R. Bravo de La Parra, L. Sanz, A discrete model of competing species sharing a parasite.}\\
\textit{ Discrete \& Continuous Dynamical Systems-B, 25(6):2121-2142, 2020.}\\
\textit{https://doi.org/10.3934/dcdsb.2019204}\\
\textit{AIMS Sciences}\\
\hline
\end{tabular}
}%

\maketitle

\centerline{\scshape RAFAEL BRAVO DE LA PARRA$^*$}
\medskip
{\footnotesize
 \centerline{U.D. Matem\'aticas, Ed. Ciencias}
   \centerline{Universidad de Alcal\'a}
   \centerline{28871 Alcal\'a de Henares, Spain}
} 

\medskip

\centerline{\scshape LUIS SANZ}
\medskip
{\footnotesize
 \centerline{Dpto. Matem\'atica Aplicada a la Ingenier\'{i}a}
   \centerline{ETSI Industriales, Univ. Polit\'ecnica de Madrid}
   \centerline{28006 Madrid, Spain}
}

\bigskip

 \centerline{(Communicated by Pierre Magal)}

\begin{abstract}
In this work we develop a discrete model of competing species affected by a
common parasite. We analyze the influence of the fast development of the shared disease on
the community dynamics.
The model is presented under the form of a two time scales discrete system with four variables. Thus, it becomes analytically tractable with the help of the appropriate reduction method.
The 2-dimensional reduced system, that has the same the asymptotic behaviour of the full model, is a generalization of the Leslie-Gower competition model. It has the unfrequent property in this kind of models of including multiple equilibrium attractors of mixed type.
The analysis of the reduced system shows that parasites can completely alter the outcome of competition depending on the parasite's basic reproductive number $R_0$. In some cases, initial conditions decide among several exclusion or coexistence scenarios.
\end{abstract}

\section{Introduction.}\label{intro}
Already in the classic experiment of Park
\cite{Park48} in the 1940s, the influence of parasites on species
competition was experimentally demonstrated. The presence of the sporozoan
parasite \emph{Adelina tribolii} can change the competition outcome of flour
beetles \emph{Tribolium confusum} and \emph{Tribolium castaneum}. In the
absence of the pathogen, \emph{Tribolium confusum} is the superior competitor
whilst in the presence of the pathogen it becomes outcompeted by
\emph{Tribolium castaneum}. The fact that parasites can alter competitive
relationships between host species implies that they can play keystone roles
in ecological communities \cite{Hatcher11}. They can alter the
outcome of an interaction, mediating coexistence or exclusion. They can also
make use of some host species as reservoirs to infect other species in which
the parasite is more virulent.

The effects of parasites on competition can be classified into those affecting
population densities and those modifying the competitive abilities of hosts.
For the first case it is used the term \emph{parasite-mediated competition}
and, for the second, \emph{parasite-modified competition}. Both are important
for host population dynamics and community structure. It is more common that
models incorporate parasite-mediated competition considering an increased
mortality as the main effect of parasitism. Nevertheless, parasite-modified
competition as an example of trait-mediated indirect interaction (TMII) also
has significant effects on community dynamics \cite{Werner03}.

A general parasitism-competition model needs to include intraspecific
competition, interspecific competition and a parasite affecting the two hosts
both in their reproduction/survival rates and in their competitive
capabilities \cite{Hatcher11}.

In \cite{Bravo17a}, two discrete time models encompassing competition
influenced by a parasite were presented. The first one corresponds to the case
of one parasite affecting demography and intraspecific competition in a single
host. The second one considers a direct competition between two different
species with one of them being the parasite's host. In this work we extend the
study to the case in which the parasite is shared by both competing species.

As it is justified in \cite{Bravo17a}, the use of timescales in this
kind of models is relevant and it is endowed with a certain novelty. It is
well known that TMIIs act on a shorter timescale than demographic effects
\cite{Bolker03,Hatcher06}. This important issue is
reflected in our model by assuming that the transitions between infected and
noninfected individuals are faster than their demographic counterparts.

In discrete time models that encompass various processes it is usual to
consider them as taking place sequentially \cite{Klepac11}. This fact
suggests a way of introducing the separation of timescales in the construction
of the model. We assume that a large number $k$ of infection-recovery episodes
take place between two demographic-competition steps. In
\cite{Sanz08,Bravo13} it is shown how to construct this
kind of discrete time models with two timescales. At the same time, a
reduction method is proposed that helps to carry out the analytical study of
the model. In this way we can circumvent the characteristic difficulty of
analysis of complete eco-epidemic competition models (see
\cite{Hatcher11} box 2.3 and references therein).

The proposed full model is represented by a 4-dimensional discrete system
whose states variables are the number of susceptible and infected individuals
in each of the two competing species. The reduction method applied to the full
model builds up a reduced 2-dimensional system approximately describing the
dynamics of the total number of individuals of both species. Then, the
asymptotic analysis of the full model can be known by studying the asymptotic
behavior of this reduced system. It turns out to be a competitive planar
discrete system \cite{Smith98} with all its non-negative solutions
tending to non-negative equilibria. A model with the same form appears in a
different context in \cite{Marva15}, where a local stability
analysis of some particular cases is performed.

In the full model, the representation of competition is based upon the
discrete Leslie-Gower model \cite{Leslie58}, that replicates in
discrete time the competition outcomes of the classical Lotka-Volterra
competition model \cite{Cushing04} through three possible scenarios: a
globally attracting coexistence equilibrium, a globally attracting exclusion
equilibrium and, finally, two attracting exclusion equilibria. The resulting
reduced system can be thought as a generalization of the Leslie-Gower model.
Nevertheless, it displays a far richer dynamics that can include multiple
equilibrium attractors of mixed type \cite{Cushing07}, i.e., at
least one coexistence equilibrium and at least one exclusion equilibrium. In
fact, there are cases whose attractors are one coexistence and one exclusion
equilibria and others with one coexistence and two exclusion equilibria. On
the other hand, some cases can be found presenting two different coexistence
equilibria as only attractors. The existence of multiple equilibrium
attractors of mixed type has been proved in structured competition models like
LPA \cite{Edmunds07} and in a non-structured Leslie-Gower competition model with Allee effects
\cite{Chow14}. The reduced model obtained in this work could
explain, through with simple nonlinearities and without introducing structure
in the populations, the observed case of multiple mixed-type attractors in
Park experiments.

The analysis of the reduced system gives insights about the possible effects
of the parasite on the competition community. It is shown that a parasite with
a large enough basic reproductive number $R_{0}$ can completely alter the
outcome of competition. It is also proved that the final situation can be
dependent on initial conditions with up to three different attracting
equilibria and their associated basins of attraction.

This work is organized as follows: In Section \ref{sec1} it is presented the
model of species competition sharing a parasite with the parasite
transmission-recovery dynamics acting faster than the demographic dynamics.
Assuming equal transmission and recovery coefficients in both species, Section
\ref{sec2} deals with the study of the model with the help of an appropriate
reduction technique. In Section \ref{sec3} different scenarios corresponding
to the parasite affecting either species growth or species competitive
abilities are discussed. The conclusions in Section \ref{conc} and the
Appendix with the proof of the mathematical results complete the manuscript.

\section{The model.}\label{sec1}
We consider two species that compete and a parasite that infects
both of them. We assume that the disease acts on a shorter timescale than the
demographic dynamics. To include two timescales in a discrete model, the time
unit of the system should be the one associated with the slow process
\cite{Bravo13,Sanz08}, in this case demography. In this
way, during one time unit we can assume that a single episode of demographic
change occurs following a number $k$ of disease infection-recovery cycles. In
this way, we intend to represent the fact that pathogens exhibit outbreaks on
short timescales, days or weeks, in comparison to demographic changes that
might be considered annual.

To represent the competition, intra and inter-species, we generalize the
well-known Leslie-Gower model \cite{Leslie58} that, for two
species $N^{1}$ and $N^{2}$, reads as follows
\begin{equation}%
\begin{array}
[c]{l}%
N^{1}(t+1)=\dfrac{b^{1}N^{1}(t)}{1+c^{11}N^{1}(t)+c^{12}N^{2}(t)}\\
\rule{0ex}{4ex}N^{2}(t+1)=\dfrac{b^{2}N^{2}(t)}{1+c^{21}N^{1}(t)+c^{22}%
N^{2}(t)},
\end{array}
\label{LGm}%
\end{equation}
where parameters $b^{i}>0$ are the inherent growth rates and $c^{ij}>0$ the
competition coefficients. In \cite{Cushing04} it is proved that this
model exhibits the same four dynamic scenarios as the classical continuous
Lotka-Volterra competition model, provided that both species have inherent
growth rates larger than one. They are characterized through the following two
parameters:
\begin{equation}
D^{1}=\frac{b^{1}-1}{c^{11}}-\frac{b^{2}-1}{c^{21}},\ D^{2}=\frac{b^{2}%
-1}{c^{22}}-\frac{b^{1}-1}{c^{12}}\label{LG-Ds}%
\end{equation}

\begin{enumerate}

\item Species 1 out-competes species 2: $D^{1}>0$ and $D^{2}<0$.

\item Species 2 out-competes species 1: $D^{1}<0$ and $D^{2}>0$.

\item Coexistence of both species: $D^{1}<0$ and $D^{2}<0$.

\item Exclusion of species either 1 or 2 depending on initial conditions:
$D^{1}>0$ and $D^{2}>0$.
\end{enumerate}

The effect of the parasite is reflected in the model by distinguishing
susceptible and infected individuals in both competing species. Let $N_{S}%
^{1}$, $N_{I}^{1}$, $N_{S}^{2}$ and $N_{I}^{2}$ denote the corresponding state
variables. The proposed extension of the Leslie-Gower model is
\begin{equation}%
\begin{array}
[c]{l}%
N_{S}^{1}(t+1)=\dfrac{b_{S}^{1}N_{S}^{1}(t)}{1+c_{SS}^{11}N_{S}^{1}%
(t)+c_{SI}^{11}N_{I}^{1}(t)+c_{SS}^{12}N_{S}^{2}(t)+c_{SI}^{12}N_{I}^{2}(t)}\\
\rule{0ex}{4ex}N_{I}^{1}(t+1)=\dfrac{b_{I}^{1}N_{I}^{1}(t)}{1+c_{IS}^{11}%
N_{S}^{1}(t)+c_{II}^{11}N_{I}^{1}(t)+c_{IS}^{12}N_{S}^{2}(t)+c_{II}^{12}%
N_{I}^{2}(t)}\\
\rule{0ex}{4ex}N_{S}^{2}(t+1)=\dfrac{b_{S}^{2}N_{S}^{2}(t)}{1+c_{SS}^{21}%
N_{S}^{1}(t)+c_{SI}^{21}N_{I}^{1}(t)+c_{SS}^{22}N_{S}^{2}(t)+c_{SI}^{22}%
N_{I}^{2}(t)}\\
\rule{0ex}{4ex}N_{I}^{2}(t+1)=\dfrac{b_{I}^{2}N_{I}^{2}(t)}{1+c_{IS}^{21}%
N_{S}^{1}(t)+c_{II}^{21}N_{I}^{1}(t)+c_{IS}^{22}N_{S}^{2}(t)+c_{II}^{22}%
N_{I}^{2}(t)},
\end{array}
\label{sd}%
\end{equation}
where all parameters are positive with $b$ representing growth rates and $c$
competition coefficients. All parameters, those representing growth as well as
those related to competition, depend on infection status. This means that we
include in the model both parasite-mediated and parasite-modified competition.

System \eqref{sd} represents the episode of demographic change. To build up
the complete model we still need to define the action of the disease
infection-recovery cycles. In order to do it we extend the discrete-time SIS
epidemic model studied in \cite{Allen94} to the case of two species
sharing the same parasite. Denoting $\mathbf{N}=\left(  N_{S}^{1},N_{I}%
^{1},N_{S}^{2},N_{I}^{2}\right)  $, the associated map is
\begin{equation}
\mathbf{F}(\mathbf{N})=\left(  F_{S}^{1}(\mathbf{N}),F_{I}^{1}(\mathbf{N}%
),F_{S}^{2}(\mathbf{N}),F_{I}^{2}(\mathbf{N})\right)  , \label{fd}%
\end{equation}
where
\[%
\begin{array}
[c]{l}%
F_{S}^{1}(\mathbf{N})=N_{S}^{1}-\dfrac{N_{S}^{1}\left(  \beta_{11}N_{I}%
^{1}+\beta_{12}N_{I}^{2}\right)  }{N_{S}^{1}+N_{I}^{1}+N_{S}^{2}+N_{I}^{2}%
}+\gamma_{1}N_{I}^{1}\\
\rule{0ex}{4ex}F_{I}^{1}(\mathbf{N})=N_{I}^{1}+\dfrac{N_{S}^{1}\left(
\beta_{11}N_{I}^{1}+\beta_{12}N_{I}^{2}\right)  }{N_{S}^{1}+N_{I}^{1}%
+N_{S}^{2}+N_{I}^{2}}-\gamma_{1}N_{I}^{1}\\
\rule{0ex}{4ex}F_{S}^{2}(\mathbf{N})=N_{S}^{2}-\dfrac{N_{S}^{2}\left(
\beta_{21}N_{I}^{1}+\beta_{22}N_{I}^{2}\right)  }{N_{S}^{1}+N_{I}^{1}%
+N_{S}^{2}+N_{I}^{2}}+\gamma_{2}N_{I}^{2}\\
\rule{0ex}{4ex}F_{I}^{2}(\mathbf{N})=N_{I}^{2}+\dfrac{N_{S}^{2}\left(
\beta_{21}N_{I}^{1}+\beta_{22}N_{I}^{2}\right)  }{N_{S}^{1}+N_{I}^{1}%
+N_{S}^{2}+N_{I}^{2}}-\gamma_{2}N_{I}^{2},
\end{array}
\]
where parameters $\beta$, representing the transmission coefficients, and
parameters $\gamma$, the recovery rates, are all positive.

As the working set we consider the possible (non-negative) population values
for which there is al least one infected individual, i.e., the set
\begin{equation}
\label{omega}\Omega:=\left\{ \left( N_{S}^{1},N_{I}^{1},N_{S}^{2},N_{I}%
^{2}\right) \in\mathbb{R}_{+}^{4}: N_{I}^{1}+N_{I}^{2}>0\right\} ,
\end{equation}
where $\mathbb{R}_{+}^{4}$ denotes the closed non-negative cone. For the
transformation to make sense in this context, the condition $\mathbf{F}%
(\Omega)\subset\Omega$ must hold. In the next section we give necessary and
sufficient conditions for this to hold in the particular case that we will
deal with.

The full model combining the demographic and the disease processes is defined
by composing the k-th iterate $\mathbf{F}^{(k)}$ of map $\mathbf{F}$
\eqref{fd}, and the map $\mathbf{S}$ associated to system \eqref{sd}:
\begin{equation}
\label{fm}\mathbf{N}(t+1)=\left(  \mathbf{S}\circ\mathbf{F}^{(k)}\right)
(\mathbf{N}(t)).
\end{equation}

To study this model we proceed in the next section to perform its reduction
into a 2-dimensional competition system with the total population of the two
species as state variables. The reduction method can be found in
\cite{Sanz08,Bravo13}. Roughly, the reduction is based on
the fact that the infectious process rapidly attains an equilibrium with fixed
proportions of susceptible and infected individuals in each species.

\section{Analysis of the reduced model with homogeneous disease transmission
and recovery.}\label{sec2}
The reduction of system \eqref{fm} in its general form is too
involved if not impossible to carry out, and so we will address a particular
case that yields interesting enough insights to be worth studying.

We assume that both species are equal in terms of disease transmission and
recovery:
\[
\beta:=\beta_{11}=\beta_{12}=\beta_{21}=\beta_{22}\text{ and }\gamma
:=\gamma_{1}=\gamma_{2}.
\]
Moreover, in order for map $\mathbf{F}$ to be properly defined and for the
disease to attaint an endemic equilibrium, we also assume the following hypothesis:

\textbf{Hypothesis 1}: $0<\gamma<\beta\leq1$.

In Appendix \ref{app1} it is proved that, under the previous hypothesis, the
disease process associated to map $\mathbf{F}$ verifies $\mathbf{F}%
(\Omega)\subset\Omega$, leaves invariant the total population of each species
$N^{1}:=N_{S}^{1}+N_{I}^{1}$ and $N^{2}:=N_{S}^{2}+N_{I}^{2}$, and attains an
equilibrium of the form:
\begin{equation}
\left(  \nu N^{1},(1-\nu)N^{1},\nu N^{2},(1-\nu)N^{2}\right)  \label{eqrap}%
\end{equation}
expressed in terms of the inverse of the parasite $R_{0}$:
\begin{equation}
\label{nu}\nu:=\frac{\gamma}{\beta}=\frac{1}{R_{0}}%
\end{equation}

The reduction procedure for system (\ref{fm}) is carried out by assuming that
this equilibrium is attained. Then it is straightforward to write the
following 2-dimensional reduced system whose state variables are the total
populations $N^{1}$ and $N^{2}$:
\begin{equation}%
\begin{array}
[c]{l}%
N^{1}(t+1)=\dfrac{r_{S}^{1}N^{1}(t)}{1+c_{S1}^{1}N^{1}(t)+c_{S2}^{1}N^{2}%
(t)}+\dfrac{r_{I}^{1}N^{1}(t)}{1+c_{I1}^{1}N^{1}(t)+c_{I2}^{1}N^{2}(t)}\\
\rule{0ex}{5ex}N^{2}(t+1)=\dfrac{r_{S}^{2}N^{2}(t)}{1+c_{S1}^{2}%
N^{1}(t)+c_{S2}^{2}N^{2}(t)}+\dfrac{r_{I}^{2}N^{2}(t)}{1+c_{I1}^{2}%
N^{1}(t)+c_{I2}^{2}N^{2}(t)},
\end{array}
\label{rm}%
\end{equation}
where for $i,j=1,2$
\[
r_{S}^{i}:=b_{S}^{i}\nu,\quad r_{I}^{i}:=b_{I}^{i}(1-\nu),\ c_{Sj}^{i}:=\nu
c_{SS}^{ij}+(1-\nu)c_{SI}^{ij},\text{ }c_{Ij}^{i}:=\nu c_{IS}^{ij}%
+(1-\nu)c_{II}^{ij}.
\]

This system generalizes the classical Leslie-Gower competition model, since it
coincides with it when $v=1$.

In Appendix \ref{app1} it is proved that the analysis of stability of the
equilibria and periodic solutions of system \eqref{fm} can be performed by
carrying out the corresponding analysis in system \eqref{rm}.

In order to carry out the mathematical treatment of system (\ref{rm}) we
express it in the form
\begin{equation}
(N^{1}(t+1),N^{2}(t+1))=H(N^{1}(t),N^{2}(t)), \label{sist}%
\end{equation}
where $H$ is the map defined by%
\[
H(x_{1},x_{2})=\left(  H_{1}(x_{1},x_{2}),H_{2}(x_{1},x_{2})\right)  =\left(
\phi_{1}(x_{1},x_{2})x_{1},\phi_{2}(x_{1},x_{2})x_{2}\right)  ,
\]
and
\[
\phi_{i}(x_{1},x_{2}):=\dfrac{r_{S}^{i}}{1+c_{S1}^{i}x_{1}+c_{S2}^{i}x_{2}%
}+\dfrac{r_{I}^{i}}{1+c_{I1}^{i}x_{1}+c_{I2}^{i}x_{2}},i=1,2.
\]
Note that all the parameters of the model are positive except possibly
$r_{I}^{1}$ and $r_{I}^{2}$ that are non-negative.

As we will see, the numbers
\[
\phi_{i}(0,0)=r_{S}^{i}+r_{I}^{i}=b_{S}^{i}\nu+b_{I}^{i}(1-\nu),\ i=1,2
\]
will play an important role in the dynamics of the system.

Let $A_{1}=\left\{  (x_{1},0):x_{1}>0\right\}  $ and $A_{2}=\left\{
(0,x_{2}):x_{2}>0\right\}  $ be the positive axes. It is immediate to realize
that the sets $\mathbf{R}_{+}^{2}$, $\mathring{\mathbf{R}}_{+}^{2}$ and
$A_{i}$, $i=1,2$ are forward invariant by $H$. In the sequel, unless otherwise
stated we will always assume that we are working on $\Omega$ \eqref{omega}.

Let us study the isoclines of the system, i.e., the sets defined by
$x_{1}=H_{1}(x_{1},x_{2})$ and $x_{2}=H_{2}(x_{1},x_{2})$. Clearly, besides
the $x_{1}$-axis (resp. $x_{2}$-axis) in which the variable $x_{2}$ (resp.
$x_{1}$) is constant, the isoclines are the curves $S_{i}$ defined by
$\phi_{i}(x_{1},x_{2})=1$ for $i=1,2$. The following Lemma presents their main properties:

\begin{lemma}
\label{prop:lema1}Let $i\in\left\{  1,2\right\}  $ be fixed. The set
$S_{i}:=\left\{  (x_{1},x_{2})\in\mathbb{R}^{2}:\phi_{i}(x_{1},x_{2}%
)=1\right\}  $ is a hyperbola that degenerates if and only if $c_{S1}%
^{i}c_{I2}^{i}=c_{S2}^{i}c_{I1}^{i}$, in which case it becomes two parallel
lines. In addition, $S_{i}$ intersects $\mathring{\mathbf{R}}_{+}^{2}$ if and
only if $\phi_{i}(0,0)>1$ and, in that case (a) only one of its branches
intersects $\mathring{\mathbf{R}}_{+}^{2}$ and (b) $S_{i}$ intersects both the
positive axes $A_{1}$ and $A_{2}$, being the intercepts defined by
\begin{equation}
R_{ij}=\frac{1}{2c_{Sj}c_{Ij}}\left(  \alpha_{j}^{i}+\sqrt{\left(  \alpha
_{j}^{i}\right)  ^{2}+4c_{Sj}^{i}c_{Ij}^{i}\left(  r_{S}^{i}+r_{I}%
^{i}-1\right)  }\right)  ,\ j=1,2, \label{intercepts}%
\end{equation}
where $\alpha_{j}^{i}:=r_{S}^{i}(c_{Ij}^{i}-1)+r_{I}^{i}(c_{Sj}^{i}%
-1),\ j=1,2$. Moreover, for each $i=1,2,$ $\Gamma_{i}:=S_{i}\cap$
$\mathbf{R}_{+}^{2}$ can be written in the form $x_{2}=\Phi_{i}(x_{1})$,
$x_{1}\in\lbrack0,R_{i1}]$, where $\Phi_{i}$ is a strictly decreasing convex function.
\end{lemma}

\begin{proof}
See Appendix \ref{app2}.
\end{proof}

In what follows we will write $\left(  x_{1},x_{2}\right)  \leq\left(
x_{1}^{\prime},x_{2}^{\prime}\right)  $ (resp. $\left(  x_{1},x_{2}\right)
<\left(  x_{1}^{\prime},x_{2}^{\prime}\right)  $) to denote that $x_{1}\leq
x_{1}^{\prime}$ and $x_{2}\leq x_{2}^{\prime}$ (resp. $x_{1}<x_{1}^{\prime}$
and $x_{2}<x_{2}^{\prime}$). Similarly, we define the $K$-order in the
following way: we write $\left(  x_{1},x_{2}\right)  \leq_{K}\left(
x_{1}^{\prime},x_{2}^{\prime}\right)  $ (resp. $\left(  x_{1},x_{2}\right)
<_{K}\left(  x_{1}^{\prime},x_{2}^{\prime}\right)  $) to denote that
$x_{1}\leq x_{1}^{\prime}$ and $x_{2}\geq x_{2}^{\prime}$ (resp. $x_{1}%
<x_{1}^{\prime}$ and $x_{2}>x_{2}^{\prime}$).

\begin{proposition}
\label{prop:prop01}Let us consider system (\ref{sist}):

a. All solutions in $\mathbf{R}_{+}^{2}$ are forward bounded, and more
specifically
\begin{equation}
H\left(  \mathbf{R}_{+}^{2}\right)  \subset S:=\left[  0,\dfrac{r_{S}^{1}%
}{c_{S1}^{1}}+\dfrac{r_{I}^{1}}{c_{I1}^{1}}\right)  \times\left[
0,\dfrac{r_{S}^{2}}{c_{S2}^{2}}+\dfrac{r_{I}^{2}}{c_{I2}^{2}}\right)
\label{eqS}%
\end{equation}

b. $H$ is strongly competitive in $\mathbf{R}_{+}^{2}$, i.e., if $x,x^{\prime
}\in\mathbf{R}_{+}^{2}$ are distinct points with $x\leq_{K}x^{\prime}$ it
follows that $H(x)<_{K}H(x^{\prime})$ \cite{Smith98}.

c. All orbits in $\mathbf{R}_{+}^{2}$ are eventually componentwise monotone,
i.e., for each \newline$\left(  N^{1}(0),N^{2}(0)\right)  \in\mathbf{R}%
_{+}^{2}$, the corresponding solution $\left(  N^{1}(t),N^{2}(t)\right)  $
verifies that $N^{i}(t)$ is eventually monotone for each $i=1,2$. Moreover,
all orbits tend to an equilibrium as $t\rightarrow\infty$.
\end{proposition}

\begin{proof}
See Appendix \ref{app2}.
\end{proof}

Let us now consider the existence of equilibriums for system (\ref{sist}).
Note that $E_{0}^{\ast}:=(0,0)$ is an equilibrium point for all values of the
parameters. Using Lemma \ref{prop:lema1} we conclude:

\begin{itemize}
\item For each $i=1,2$, there exists a semitrivial equilibrium point
$E_{i}^{\ast}$ on the positive axis $A_{i}$ if and only if $\phi_{i}(0,0)>1$.
In that case the semitrivial equilibrium is unique.

\item A necessary condition for the existence of a positive equilibrium is
that $\phi_{i}(0,0)>1$ for both $i=1,2$.
\end{itemize}

The next result analyzes the global behavior of solutions of system
(\ref{sist}) except in the case in which $\phi_{i}(0,0)>1$, $i=1,2$:

\begin{theorem}
\label{prop:prop2}Let us consider system (\ref{sist}):

a. For each $i=1,2$, if $\phi_{i}(0,0)\leq1$ then for any initial value on
$\mathbf{R}_{+}^{2}$ species $i$ tends to extinction as $t\rightarrow\infty$.

b. If $\phi_{i}(0,0)\leq1$ for $i=1,2$, all orbits in $\mathbf{R}_{+}^{2}$
tend to $E_{0}^{\ast},$ and if the inequalities are strict then $E_{0}^{\ast}
$ is hyperbolic. If $\phi_{i}(0,0)>1$ for $i=1,2,$ then $E_{0}^{\ast}$ is a
repeller and therefore no orbit can converge to $E_{0}^{\ast}$.

c. If $\phi_{1}(0,0)>1$ and $\phi_{2}(0,0)\leq1$ then (1) all orbits with
$N^{1}(0)=0$ tend to $E_{0}^{\ast}$ and (2) all orbits with $N^{1}(0)>0$ tend
to $E_{1}^{\ast}$.

d. If $\phi_{2}(0,0)>1$ and $\phi_{1}(0,0)\leq1$ then (1) all orbits with
$N^{2}(0)=0$ tend to $E_{0}^{\ast}$ and (2) all orbits with $N^{2}(0)>0$ tend
to $E_{2}^{\ast}$.
\end{theorem}

\begin{proof}
See Appendix \ref{app2}.
\end{proof}

\smallskip

Let us now consider the case in which $\phi_{i}(0,0)>1$ for $i=1,2$. The
positive equilibria of the system are the (positive) number of intersections
of $\Gamma_{1}$ and $\Gamma_{2}$. In the first place, since the isoclines are
hyperbolas they can intersect in four points at most, and so there can be at
most four positive equilibria. It is easy to find parameter values for which
there are zero, one, two or three positive equilibria. After extensive
numerical simulations we have not been able to find any case for which there
are four positive intersections, and so we will not treat that case in the
following discussion.

In order to study the behavior of solutions, we consider the different generic
cases based on the relative position of the intercepts $R_{ij}$
(\ref{intercepts}) of $\Gamma_{1}$ and $\Gamma_{2},$ and on the number (up to
three) of positive equilibriums. Specifically, we distinguish several
scenarios that we denote with a letter that corresponds to the relative
position of the intercepts $R_{ij}$ and a subindex that denotes the number of
positive equilibriums. We only consider the generic cases, i.e., we omit the
cases in which the isoclines are tangent at an equilibrium point:%
\begin{equation}%
\begin{tabular}
[c]{ll}%
- \textbf{Case A}. $R_{11}<R_{21},$ $R_{12}>R_{22}$. & There can be one
(\textbf{Case A}$_{1}$)\\
& or three (\textbf{Case A}$_{3}$) positive equilibria.\\
- \textbf{Case B}. $R_{11}>R_{21},$ $R_{12}<R_{22}$. & There can be one
(\textbf{Case B}$_{1}$)\\
& or three (\textbf{Case B}$_{3}$) positive equilibria.\\
- \textbf{Case C}. $R_{11}>R_{21},$ $R_{12}>R_{22}$. & There can be zero
(\textbf{Case C}$_{0}$)\\
& or two (\textbf{Case C}$_{2}$) positive equilibria.\\
- \textbf{Case D}. $R_{11}<R_{21},$ $R_{12}<R_{22}$. & There can be zero
(\textbf{Case D}$_{0}$)\\
& or two (\textbf{Case D}$_{2}$) positive equilibria.
\end{tabular}
\ \label{casos}%
\end{equation}

Taking into account Lemma \ref{prop:lema1}, the isoclines divide
$\mathring{\mathbf{R}}_{+}^{2}$ in a finite number of open connected and
disjoint sets in which $x_{i}$ is either strictly increasing or strictly
decreasing for each $i=1,2$. To describe the kind of monotonicity in each
region we use arrows, in such a way that, for example, the situation in which
$x_{1}$ decreases and $x_{2}$ increases corresponds to a horizontal arrow
pointing to the left and a vertical arrow pointing up.

\begin{figure}[htp]
\centering
\begin{subfigure}[t]{0.32\textwidth}
\centering
\includegraphics[width=\linewidth]{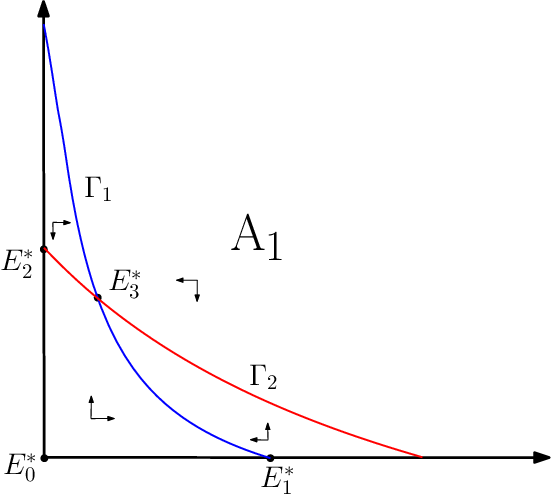}
\end{subfigure}
\hfill\begin{subfigure}[t]{0.33\textwidth}
\centering
\includegraphics[width=\linewidth]{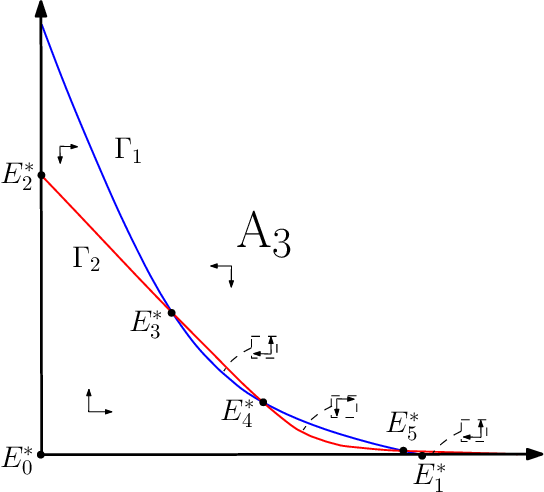}
\end{subfigure}
\hfill\begin{subfigure}[t]{0.32\textwidth}
\centering
\includegraphics[width=\linewidth]{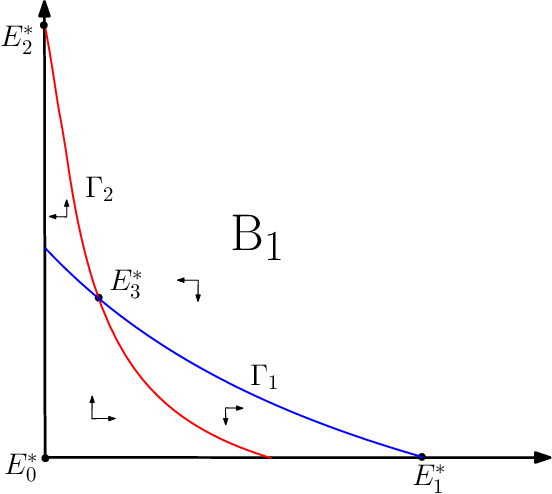}
\end{subfigure}
\par
\rule{0ex}{2ex}
\par
\centering
\begin{subfigure}[t]{0.33\textwidth}
\centering
\includegraphics[width=\linewidth]{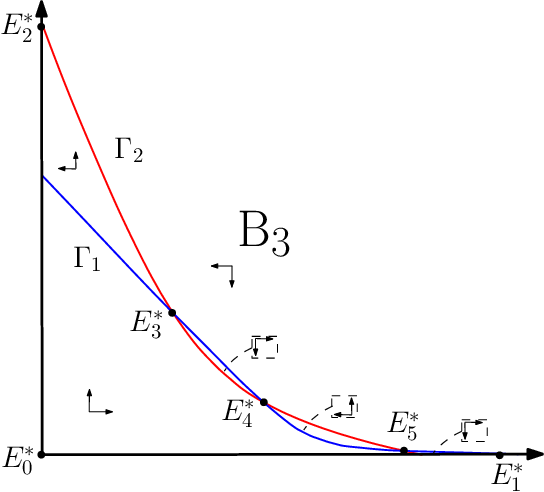}
\end{subfigure}
\hfill\begin{subfigure}[t]{0.32\textwidth}
\centering
\includegraphics[width=\linewidth]{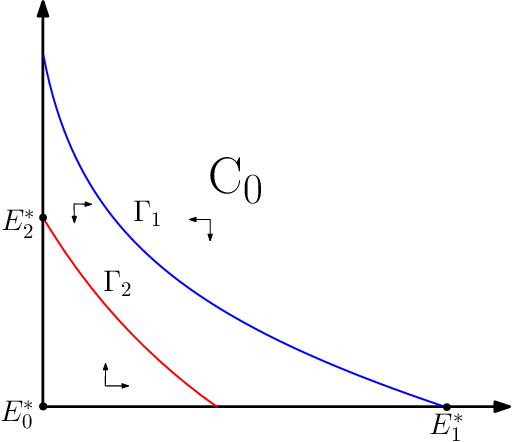}
\end{subfigure}
\hfill\begin{subfigure}[t]{0.32\textwidth}
\centering
\includegraphics[width=\linewidth]{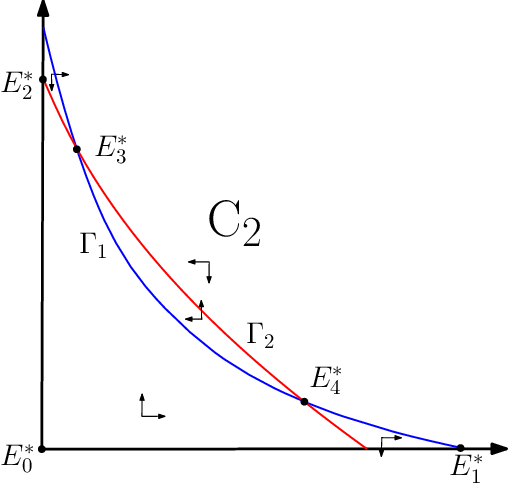}
\end{subfigure}
\par
\rule{0ex}{2ex}
\par
\centering
\begin{subfigure}[t]{0.33\textwidth}
\centering
\includegraphics[width=\linewidth]{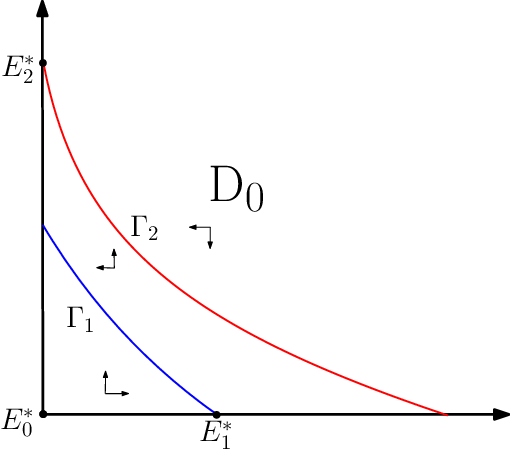}
\end{subfigure}
\qquad\begin{subfigure}[t]{0.33\textwidth}
\centering
\includegraphics[width=\linewidth]{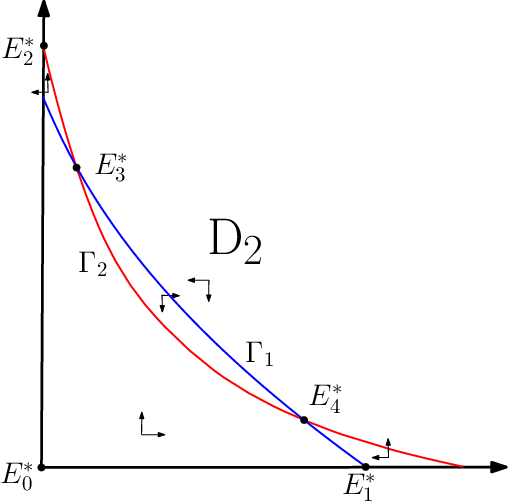}
\end{subfigure}
\caption{Different configurations of system (\ref{sist}) when $\phi
_{i}(0,0)>1$ for $i=1,2$, in terms of the relative position of the intercepts
of isoclines, $R_{ij}$ \eqref{intercepts} and the number of positive
equilibria, as described in (\ref{casos})}
\label{fig1}
\end{figure}

Figure \ref{fig1} shows the different configurations described in (\ref{casos}).

The following result deals on the one hand with the stability of the
semitrivial and of the positive equilibria in the different scenarios
described in (\ref{casos}). On the other hand, it characterizes the global
dynamics on the positive axes and studies the behavior of solutions for
positive initial conditions:

\begin{theorem}
\label{prop:prop3}Let us consider system (\ref{sist}) and assume that
$\phi_{i}(0,0)>1$ for $i=1,2$.

\textbf{1. Semitrivial equilibria and behavior on the positive axes.}

\qquad a. All orbits starting on the positive $i$-axis $A_{i}$ converge to
$E_{i}^{\ast}$, $i=1,2$.

\qquad b. Let $i=1,2$ be fixed. If $\phi_{j}(E_{i}^{\ast})<1$ for $j\neq i,$
then $E_{i}^{\ast}$ is hyperbolic and attracting. If $\phi_{j}(E_{i}^{\ast
})>1$ for $j\neq i$ then $E_{i}^{\ast}$ is hyperbolic and unstable.

\textbf{2. Position of the equilibria.} Let us assume that there exist
$s\geq0$ positive equilibria that we denote $E_{2+i}^{\ast}$, $i=1,...,s$.
Then they can be ordered using the $K$-order in such a way that%
\[
E_{2}^{\ast}<_{K}E_{3}^{\ast}<_{K}\cdots<_{K}E_{2+s}^{\ast}<_{K}E_{1}^{\ast}%
\]
and in the sequel we assume that they are ordered in this way.

\textbf{3. Stability of equilibria and behavior of solutions for positive
initial conditions. }Let us consider cases \textbf{A}, \textbf{B}, \textbf{C}
and \textbf{D} as described in (\ref{casos}):

\begin{itemize}
\item \textbf{Case} \textbf{A}. $E_{1}^{\ast}$ and $E_{2}^{\ast}$ are
hyperbolic and unstable.

\begin{itemize}
\item In case \textbf{A}$_{1},$ the (unique) positive equilibrium $E_{3}%
^{\ast}$ is hyperbolic and attracting. All orbits starting on $\mathring
{\mathbf{R}}_{+}^{2}$ converge to $E_{3}^{\ast}$.

\item In case \textbf{A}$_{3}$, $E_{3}^{\ast}$ and $E_{5}^{\ast}$ are
hyperbolic and attracting and $E_{4}^{\ast}$ is a hyperbolic saddle. Moreover
$E_{4}^{\ast}$ can not attract any open set and therefore almost all orbits
starting in $\mathring{\mathbf{R}}_{+}^{2}$ converge to either $E_{3}^{\ast}$
or $E_{5}^{\ast}$.
\end{itemize}

\item \textbf{Case B}. $E_{1}^{\ast}$ and $E_{2}^{\ast}$ are hyperbolic and attracting.

\begin{itemize}
\item In case \textbf{B}$_{1},$ the (unique) positive equilibrium $E_{3}%
^{\ast}$ is a hyperbolic saddle. Moreover $E_{3}^{\ast}$ can not attract any
open set and therefore almost all orbits starting in $\mathring{\mathbf{R}%
}_{+}^{2}$ converge to either $E_{1}^{\ast}$ or $E_{2}^{\ast}$.

\item In case \textbf{B}$_{3},$ $E_{3}^{\ast}$ and $E_{5}^{\ast}$ are
hyperbolic saddles and $E_{4}^{\ast}$ is hyperbolic and attracting. Moreover
$E_{3}^{\ast}$ and $E_{5}^{\ast}$ can not attract any open set and therefore
almost all orbits starting in $\mathring{\mathbf{R}}_{+}^{2}$ converge to
either $E_{1}^{\ast}$, $E_{2}^{\ast}$ or $E_{4}^{\ast}$.
\end{itemize}

\item \textbf{Case C}. $E_{1}^{\ast}$ is hyperbolic and attracting and
$E_{2}^{\ast}$ is hyperbolic and unstable.

\begin{itemize}
\item In case \textbf{C}$_{0}$ all orbits starting on $\mathring{\mathbf{R}%
}_{+}^{2}$ converge to $E_{1}^{\ast}$.

\item In case \textbf{C}$_{2}$ no orbits starting in $\mathring{\mathbf{R}%
}_{+}^{2}$ can converge to $E_{2}^{\ast}$. Besides $E_{3}^{\ast}$ is
hyperbolic and attracting and $E_{4}^{\ast}$ is hyperbolic and a saddle.
Moreover $E_{4}^{\ast}$ can not attract any open set and therefore almost all
orbits starting in $\mathring{\mathbf{R}}_{+}^{2}$ converge to either
$E_{1}^{\ast}$ or $E_{3}^{\ast}$.
\end{itemize}

\item \textbf{Case D}. $E_{2}^{\ast}$ is hyperbolic and attracting and
$E_{1}^{\ast}$ is hyperbolic and unstable.

\begin{itemize}
\item In case \textbf{D}$_{0}$ all orbits starting on $\mathring{\mathbf{R}%
}_{+}^{2}$ converge to $E_{2}^{\ast}$.

\item In case \textbf{D}$_{2}$ no orbits starting in $\mathring{\mathbf{R}%
}_{+}^{2}$ can converge to $E_{1}^{\ast}$. Besides $E_{4}^{\ast}$ is
hyperbolic and attracting and $E_{3}^{\ast}$ is hyperbolic and a saddle.
Moreover $E_{3}^{\ast}$ can not attract any open set and therefore almost all
orbits starting in $\mathring{\mathbf{R}}_{+}^{2}$ converge to either
$E_{2}^{\ast}$ or $E_{4}^{\ast}$.
\end{itemize}
\end{itemize}
\end{theorem}

\begin{proof}
See Appendix \ref{app2}.
\end{proof}

\section{Parasite-mediated and parasite-modified competition.}

\label{sec3} In this section we discuss the effects of parasitism on the
competition community, system \eqref{fm}, with the help of the analysis of
system \eqref{rm} performed in Section \ref{sec2}. The results on approximate
aggregation in \cite{Bravo13,Sanz08} guarantee, loosely
speaking, that when a solution of the reduced system (\ref{rm}) with initial
condition $(N^{1}(0),N^{2}(0))$ tends to a hyperbolic equilibrium $(N^{1\ast
},N^{2\ast})$ then, for $k$ large enough, any positive solution of system
(\ref{fm}) verifying $N_{S}^{1}(0)+N_{I}^{1}(0)=N^{1}(0)$ and $N_{S}%
^{2}(0)+N_{I}^{2}(0)=N^{2}(0)$ tends to an equilibrium which is approximately
\[
\left(  \nu N^{1\ast},\left(  1-\nu\right)  N^{1\ast},\nu N^{2\ast},\left(
1-\nu\right)  N^{2\ast}\right)  ,
\]
with $\nu$ given by (\ref{nu}).

We recall that the Leslie-Gower competition model presents the same four
dynamic scenarios as the continuous Lotka-Volterra competition model: 1.
Species 1 out-competes species 2; 2. Species 2 out-competes species 1; 3.
Coexistence; 4. Exclusion of either species 1 or 2. These four scenarios are
also among the options described in Theorem \ref{prop:prop3} for system
(\ref{rm}). Case C$_{0}$ corresponds to 1., case D$_{0}$ to 2., case A$_{1}$
to 3. and case B$_{1}$ to 4..

\begin{figure}[htp]
\begin{center}
\includegraphics[width=0.5\linewidth]{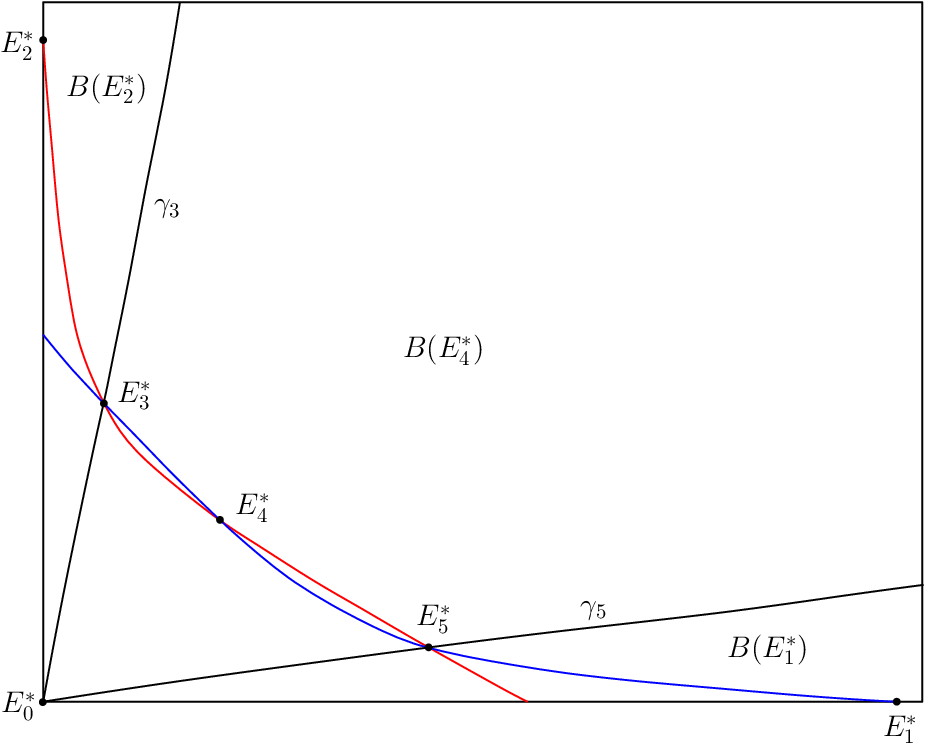}
\caption{Basins of attraction $B(E_{1}^{\ast})$, $B(E_{2}^{\ast})$ and $B(E_{4}^{\ast})$ of
equilibria $E_{1}^{\ast}$, $E_{2}^{\ast}$ and $E_{4}^{\ast}$ and separatrix
curves $\gamma_{3}$ and $\gamma_{5}$ for system \eqref{rm} for parameters
values: $\nu=0.5$, $b_{S}^{1}=13$, $b_{I}^{1}=3.6$, $b_{S}^{2}=3.4$,
$b_{I}^{2}=8$, $c_{SS}^{11}=c_{SI}^{11}=0.9$, $c_{IS}^{11}=c_{II}^{11}=0.1$,
$c_{SS}^{12}=c_{SI}^{12}=1.1$, $c_{IS}^{12}=c_{II}^{12}=5$, $c_{SS}%
^{21}=c_{SI}^{21}=6$, $c_{IS}^{21}=c_{II}^{21}=0.3$, $c_{SS}^{22}=c_{SI}%
^{22}=0.2$, $c_{IS}^{22}=c_{II}^{22}=0.8$.}
\label{figCA}
\end{center}
\end{figure}

The existence of a parasite affecting one of the competing species
\cite{Bravo17a} makes a new scenario appear consisting in either
exclusion of the uninfected species or species coexistence, i.e., parasite
mediated coexistence. The cases C$_{2}$ and D$_{2}$ in Theorem
\ref{prop:prop3} also exhibit this dichotomy between exclusion of one of the
species and coexistence depending on initial conditions. The remaining two
cases in Theorem \ref{prop:prop3}, A$_{3}$ and B$_{3}$, are specific
consequences of the existence of a parasite shared by two competing species.
They add two new dynamic scenarios due to the existence of three positive
equilibria. In case A$_{3}$, two of them are asymptotically stable whereas the
intermediate one is a saddle, and so for most initial conditions the long term
output of competition is one of those two different levels of stable
coexistence. Finally, in case B$_{3}$, it is the intermediate equilibrium that
is asymptotically stable together with the two exclusion equilibria. We note
that in this latter case all three possibilities of exclusion can be attained.
It suffices to start dynamics in the corresponding basin of attraction, see
Figure \ref{figCA}, to get species 2 excluded ($B(E_{1}^{\ast})$), to get
species 1 excluded ($B(E_{2}^{\ast})$) or to obtain coexistence ($B(E_{4}%
^{\ast})$).

The effects of parasite on competition are classified
\cite{Hatcher11} into parasite-mediated competition and
parasite-modified competition. The first one refers to parasites influencing
the competitive relationship by altering host densities through changes in
either parasite-induced mortality or fecundity. On the other hand, the term
parasite-modified competition is used for the case of parasites altering the
competitive abilities of individuals of both species.

To illustrate the parasite-mediated competition we assume that parasites do
not affect competitive abilities. Competition coefficients do not depend on
individuals being infected or uninfected. To be precise, we consider in system
\eqref{fm} four different competition coefficients just depending on the
species involved: $c^{11}$, $c^{12}$, $c^{21}$ and $c^{22}$. The sixteen
competition coefficients appearing in \eqref{fm} verify:
\[
c_{SS}^{ij}=c_{SI}^{ij}=c_{IS}^{ij}=c_{II}^{ij}=c^{ij}\ ,\text{ for }%
i,j\in\{1,2\}.
\]

Concerning the growth rates we generally assume that they are reduced by
parasites:
\[
b_{I}^{i}<b_{S}^{i}\ ,\text{ for }i\in\{1,2\}.
\]

The competition model without parasites has the following Leslie-Gower form
\eqref{LGm}:
\begin{equation}
N_{S}^{1}(t+1)=\dfrac{b_{S}^{1}N_{S}^{1}(t)}{1+c^{11}N_{S}^{1}(t)+c^{12}%
N_{S}^{2}(t)}\ ,\ N_{S}^{2}(t+1)=\dfrac{b_{S}^{2}N_{S}^{2}(t)}{1+c^{21}%
N_{S}^{1}(t)+c^{22}N_{S}^{2}(t)}.\label{pmcM1}%
\end{equation}
The reduced model \eqref{rm} also has a Leslie-Gower form:
\begin{equation}
N^{1}(t+1)=\dfrac{\left(  \nu b_{S}^{1}+(1-\nu)b_{I}^{1}\right)  N^{1}%
(t)}{1+c^{11}N^{1}(t)+c^{12}N^{2}(t)}\ ,\ N^{2}(t+1)=\dfrac{\left(  \nu
b_{S}^{2}+(1-\nu)b_{I}^{2}\right)  N^{2}(t)}{1+c^{21}N^{1}(t)+c^{22}N^{2}%
(t)},\label{pmcM2}%
\end{equation}
the only change concerning the growth rates. In the uninfected community they
are $b_{S}^{1}$ and $b_{S}^{2}$ and when the parasite is taken into account
they become
\[
b^{1}:=\nu b_{S}^{1}+(1-\nu)b_{I}^{1}\ \text{ and }\ b^{2}:=\nu b_{S}%
^{2}+(1-\nu)b_{I}^{2}.
\]

An immediate conclusion is that the parasite can drive extinct any of the
species. Indeed, a necessary condition for the species $i$ to get extinct by
the effect of the parasite is $b_{I}^{i}<1$. If this is the case, a large
enough $R_{0}$ yields
\[
b^{i}=\nu b_{S}^{i}+(1-\nu)b_{I}^{i}=\frac{1}{R_{0}}b_{S}^{i}+(1-\frac
{1}{R_{0}})b_{I}^{i}<1,
\]
and then Theorem \ref{prop:prop2} implies the long-term extinction of species
$i$.

We next assume that, independently of the size of $R_{0}$, the parasite cannot
drive extinct any of the two species, i.e., $b_{I}^{1},b_{I}^{2}>1$. The
outcome of competition in the Leslie-Gower system \eqref{pmcM2} depends (see
\eqref{LG-Ds}) on the signs of coefficients%
\[%
\begin{array}
[c]{l}%
\bar{D}^{1}=\dfrac{\nu b_{S}^{1}+(1-\nu)b_{I}^{1}-1}{c^{11}}-\dfrac{\nu
b_{S}^{2}+(1-\nu)b_{I}^{2}-1}{c^{21}},\\
\rule{0ex}{4ex}\bar{D}^{2}=\dfrac{\nu b_{S}^{2}+(1-\nu)b_{I}^{2}-1}{c^{22}%
}-\dfrac{\nu b_{S}^{1}+(1-\nu)b_{I}^{1}-1}{c^{12}},
\end{array}
\]
that can be expressed in the following form:
\[%
\begin{array}
[c]{l}%
\bar{D}^{1}=(1-\dfrac{1}{R_{0}})\left(  \dfrac{b_{I}^{1}-1}{c^{11}}%
-\dfrac{b_{I}^{2}-1}{c^{21}}\right)  +\dfrac{1}{R_{0}}\left(  \dfrac{b_{S}%
^{1}-1}{c^{11}}-\dfrac{b_{S}^{2}-1}{c^{21}}\right)  ,\\
\rule{0ex}{4ex}\bar{D}^{2}=(1-\dfrac{1}{R_{0}})\left(  \dfrac{b_{I}^{2}%
-1}{c^{22}}-\dfrac{b_{I}^{1}-1}{c^{12}}\right)  +\dfrac{1}{R_{0}}\left(
\dfrac{b_{S}^{2}-1}{c^{22}}-\dfrac{b_{S}^{1}-1}{c^{12}}\right)  .
\end{array}
\]
A straightforward conclusion of these expressions is that, for large enough
$R_{0},$ their signs coincide with those of parameters
\[
\bar{D}_{I}^{1}:=\dfrac{b_{I}^{1}-1}{c^{11}}-\dfrac{b_{I}^{2}-1}{c^{21}%
}\ ,\ \bar{D}_{I}^{2}:=\dfrac{b_{I}^{2}-1}{c^{22}}-\dfrac{b_{I}^{1}-1}{c^{12}%
},
\]
and are independent of those of parameters
\[
\bar{D}_{S}^{1}:=\dfrac{b_{S}^{1}-1}{c^{11}}-\dfrac{b_{S}^{2}-1}{c^{21}%
}\ ,\ \bar{D}_{S}^{2}:=\dfrac{b_{S}^{2}-1}{c^{22}}-\dfrac{b_{S}^{1}-1}{c^{12}%
}.
\]
Therefore, an endemic parasite establishing a large enough fraction of
infected individuals in the population can modify any outcome of competition
and yield any other. Indeed, the outcome of system \eqref{pmcM1} plays no role
in the outcome of system \eqref{pmcM2}, which would coincide with that of the
following system
\[
N^{1}(t+1)=\dfrac{b_{I}^{1}N^{1}(t)}{1+c^{11}N^{1}(t)+c^{12}N^{2}%
(t)}\ ,\ N^{2}(t+1)=\dfrac{b_{I}^{2}N^{2}(t)}{1+c^{21}N^{1}(t)+c^{22}N^{2}%
(t)}.
\]

\begin{figure}[htp]
\begin{center}
\includegraphics[scale=0.5]{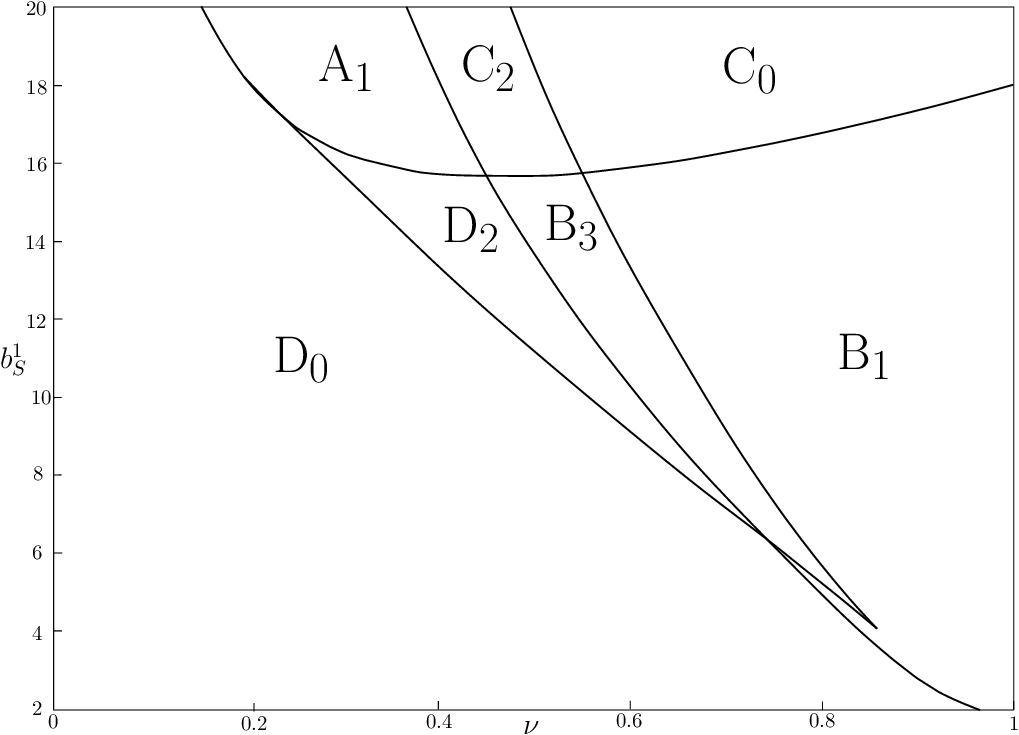}
\caption{Asymptotic behaviour cases of solutions of system (\ref{rm}) (Th. \eqref{prop:prop3}) for parameters values: $\nu\in(0,1)$, $b_{S}^{1}\in[2,20]$, $b_{I}^{1}=2$,
$b_{S}^{2}=4.4,b_{I}^{2}=9$, $c_{SS}^{11}=1.3$, $c_{SI}^{11}=0.5$,
$c_{IS}^{11}=c_{II}^{11}=0.1$, $c_{SS}^{12}=1$, $c_{SI}^{12}=0.05$,
$c_{IS}^{12}=8$, $c_{II}^{12}=3$, $c_{SS}^{21}=6$, $c_{SI}^{21}=c_{IS}%
^{21}=c_{II}^{21}=0.3$, $c_{SS}^{22}=c_{SI}^{22}=0.2$, $c_{IS}^{22}%
=c_{II}^{22}=0.8$.}%
\label{fig:bifurc}%
\end{center}
\end{figure}

Let us now illustrate how parasite-modified competition can yield a rich
variety of different outcomes. For that we develop a particular case of system
\eqref{fm}. In Figure \ref{fig:bifurc} we show the cases described in Theorem
\eqref{prop:prop3} that correspond to different values of parameter
$\nu=1/R_{0}\in(0,1)$ and parameter $b_{S}^{1}\in\lbrack2,20]$ setting the
following fixed values for the rest of parameters: $b_{I}^{1}=2$, $b_{S}%
^{2}=4.4,b_{I}^{2}=9$, $c_{SS}^{11}=1.3$, $c_{SI}^{11}=0.5$, $c_{IS}%
^{11}=c_{II}^{11}=0.1$, $c_{SS}^{12}=1$, $c_{SI}^{12}=0.05$, $c_{IS}^{12}=8$,
$c_{II}^{12}=3$, $c_{SS}^{21}=6$, $c_{SI}^{21}=c_{IS}^{21}=c_{II}^{21}=0.3$,
$c_{SS}^{22}=c_{SI}^{22}=0.2$, $c_{IS}^{22}=c_{II}^{22}=0.8$. We have chosen
the parameters values so that the competitive abilities of infected
individuals are less than or equal to those of uninfected individuals in the
same circumstances, i.e., $c_{AI}^{ij}\leq c_{AS}^{ij}$ for $A\in\{S,I\}$ and
$i,j\in\{1,2\}$.

The corresponding reduced system \eqref{rm} is

\begin{equation}
\begin{array}{l}
 \begin{split}
N^{1}(t+1)  = & \dfrac{ b_{S}^{1}\nu N^{1}(t)} {1+\left( 1.3\nu
+0.5(1-\nu)\right) N^{1}(t)+\left( \nu+0.05(1-\nu)\right) N^{2}(t)}\\
  & \quad +\dfrac{2(1-\nu)N^{1}(t)} {1+\left(
0.1\nu+0.1(1-\nu)\right) N^{1}(t)+\left( 8\nu+3(1-\nu)\right) N^{2}(t)},
 \end{split}
\\ \rule{0ex}{7ex}
 \begin{split}
 N^{2}(t+1)  = & \dfrac{ 4.4\nu N^{1}(t)} {1+\left(
6\nu+0.3(1-\nu)\right) N^{1}(t)+\left( 0.2\nu+0.2(1-\nu)\right) N^{2}(t)}\\
 &   +\dfrac{9(1-\nu)N^{1}(t)} {1+\left(
0.3\nu+0.3(1-\nu)\right) N^{1}(t)+\left( 0.8\nu+0.8(1-\nu)\right) N^{2}(t)}.
 \end{split}
\end{array}
\label{rmej}
\end{equation}

When there is no parasite (in the limit when $\nu$ tends to $1$) we find two
different scenarios. If $b_{S}^{1}\in\lbrack2,18)$ the corresponding case is
B$_{1}$ (see Figure \ref{fig1} to follow the cases) that entails an exclusion
situation of one of the species depending on initial conditions. For
$b_{S}^{1}>18$ the new scenario is case C$_{0}$ in which species 1
out-competes species 2. The influence of parasites can now be followed by
moving backward from $\nu$ close to $1$ ($R_{0}=1$) to $\nu=0$ ($R_{0}%
\rightarrow\infty$) in the bifurcation diagram (Figure \ref{fig:bifurc}).

Low values of $b_{S}^{1}$ imply, as $R_{0}$ increases, moving from case
B$_{1}$, exclusion, to case D$_{0}$ in which species 2 out-competes species 1.
We see that in this example the parasite favours the second species.

The intermediate values of $b_{S}^{1}$ yield the same beginning and end of the
dynamic scenario, B$_{1}$ and D$_{0}$, though going first through case B$_{3}$
and then case D$_{2}$. In case B$_{3}$, apart from the two exclusion options
included in case B$_{1}$, it appears the possibility, all three depending on
initial conditions (see Figure \ref{fig1}), of coexistence. This is a neat
example of parasite mediated coexistence. The transition from case B$_{3}$ to
D$_{0}$ is done through D$_{2}$ where the option of species 1 excluding
species 2 has disappeared.

High values of $b_{S}^{1}$, the growth rate of non-infected individuals of
species 1, entail a transition from case C$_{0}$, species 1 out-competes
species 2, to D$_{0}$, species 2 out-competes species 1, as $R_{0}$ increases.
This transition is done through case C$_{2}$, that encompasses species 2
exclusion and coexistence, followed by case A$_{1}$, that represent global coexistence.

These transitions from an scenario to another one as $\nu$ decreases from one
to zero can even present a more complex pattern. For example, for a value of
$b_{S}^{1}$ around $16$, the dynamic scenario goes through the following six
cases in order: B$_{1}$, C$_{0}$, C$_{2}$, A$_{1}$, D$_{2}$ and D$_{0}$.
Similarly, for values of $b_{S}^{1}$ close to $15.7$, the sequence of cases
is: B$_{1}$, B$_{3}$, C$_{2}$, A$_{1}$, D$_{2}$ and D$_{0}$.

\section{Conclusion.}

\label{conc} In this work we have extended a previous eco-epidemic model
\cite{Bravo17a} that filled a gap in the existent literature due to
being in discrete time and having the parasite dynamics occurring on a shorter
time scale than the competition interactions. In the model studied here we
have introduced the three ingredients demanded by general parasitism-competition modules \cite{Hatcher11}:
intraspecific host competition, interspecific competition between hosts and a
parasite capable of infecting both hosts. Whereas in \cite{Bravo17a}
the parasite affected only one the two competing species, the contribution of
this work is analyzing the effect of a shared parasite affecting the
competition community.

After the general presentation of the model in Section \ref{sec1}, we have
assumed in Section \ref{sec2} homogeneous disease transmission and recovery to
keep this work simple enough to be analytically tractable. In this way we have
been able to apply a reduction technique that allows one to find analytic
results stressing the effects of the disease on the community. In spite of
this simplifying assumption, the mathematical analysis of the model reveals a
number of interesting long term scenarios that were not exhibited in previous
models \cite{Bravo17a}. In particular, there are two new situations
in which three different positive equilibria exist. They are called, see
Figure \ref{fig1}, cases A$_{3}$ and B$_{3}$ in Section \ref{sec2}.

Case A$_{3}$ is a notable example of coexistence of non-Lotka-Volterra type.
Two positive equilibria out of three are asymptotically stable, whereas the
third, the intermediate one, is a saddle whose stable manifold serves as
separatrix curve of the basins of attraction of the first two. The two
exclusion equilibria are also saddles. Thus, we have that different positive
initial conditions lead to one of the two different stable coexistence equilibria.

Case B$_{3}$ is another example of coexistence of non-Lotka-Volterra type. In
this case, the intermediate equilibrium is asymptotically stable and the two
other positive equilibria are saddles (see Figure \ref{figCA}). The two
exclusion equilibria are asymptotically stable. The stable manifolds of the
saddles are the boundary of the three basins of attraction included in the
positive quadrant. Now different positive initial conditions lead to one of
the three exclusion possibilities: species 1 is the superior competitor;
species 2 is the superior competitor; species coexist. In
\cite{Leslie68} it is described a competition experiment in which
varying the initial conditions two \emph{Tribolium} species could coexist
though in most cases exclusion was the outcome. The attempt to interpret these
results in terms of the classical Lotka-Volterra scenario, or rather of its
discrete counterpart, the Leslie-Gower model, failed. These results were
reviewed in \cite{Edmunds03} and an explanation of both species
exclusion and species coexistence provided by means of discrete, non-linear
stage-structured model for flour beetles \cite{Cushing02} extended to
include two interacting species \cite{Edmunds01}. The case B$_{3}$ of
the simple planar system \eqref{rm} offers an alternative explanation in a
non-structured setting, of the simplified action of a shared parasite.

The results in Section \ref{sec3} show that parasites can reverse the effects
of direct competition between the hosts. System \eqref{pmcM1}, competition
without parasite, can yield one of the species as superior competitor whereas
the outcome of system \eqref{pmcM2}, competition with parasite only affecting
growing terms, can be the exclusion of this same species provided that the
$R_{0}$ of the parasite is large enough. In the same direction, it is
illustrated in the example corresponding to system \eqref{rmej} that the
species 2, favoured by the parasite, can overcome the growing advantage of the
species 1 as the parasite $R_{0}$ increases. From a community ecology
perspective, these results confirm the extended opinion
\cite{Hatcher11} that parasites can act as keystone species. A
preferential application of this fact is to biological control. The insights
gained with the analysis of this kind of models can lead to the fruitful use
of parasites in biological control. A more detailed analysis of the model in
particular cases relevant to specific situations and data sets, could be a
sound perspective of this work.

\bigskip

\appendix

\textbf{{\Large Appendix}}

\section{Disease dynamics in the case $\beta_{ij}=\beta$ and $\gamma
_{i}=\gamma$ for $i,j=1,2$, and reduction results.}

\label{app1} The disease dynamics keeps constant the total population of both
species
\[
N^{1}:=F_{S}^{1}(\mathbf{N})+F_{I}^{1}(\mathbf{N})=N_{S}^{1}+N_{I}^{1}%
,\ N^{2}:=F_{S}^{2}(\mathbf{N})+F_{I}^{2}(\mathbf{N})=N_{S}^{2}+N_{I}^{2}.
\]
Thus
\[
N_{S}^{1}=N^{1}-N_{I}^{1}\text{ and }\ N_{S}^{2}=N^{2}-N_{I}^{2},
\]
and the disease dynamics, for fixed $N^{1}$ and $N^{2}$, is completely defined
in terms of variables $N_{I}^{1}$ and $N_{I}^{2}$. The maps updating the
values of these variables in the particular case of equal transmission and
recovery parameters are the following :
\[
\begin{array}{l}
\begin{split}
f(N_{I}^{1},N_{I}^{2}):= & F_{I}^{1}(N^{1}-N_{I}^{1},N_{I}^{1},N^{2}-N_{I}
^{2},N_{I}^{2}) \\
 & \qquad \qquad =(1-\gamma)N_{I}^{1}+\dfrac{\beta(N^{1}-N_{I}^{1})\left(
N_{I}^{1}+N_{I}^{2}\right)}{N^{1}+N^{2}}
\end{split}\\
\begin{split}
g(N_{I}^{1},N_{I}^{2}):= & F_{I}^{2}(N^{1}-N_{I}^{1},N_{I}
^{1},N^{2}-N_{I}^{2},N_{I}^{2})\\
 & \qquad \qquad =(1-\gamma)N_{I}^{2}+\dfrac{\beta(N^{2}
-N_{I}^{2})\left(  N_{I}^{1}+N_{I}^{2}\right)  }{N^{1}+N^{2}}.
\end{split}
\end{array}
\]
As we noted when introducing the disease dynamics, we must impose some
conditions on parameters $\beta$ and $\gamma$ so that $\mathbf{F}%
(\Omega)\subset\Omega$ \eqref{omega}. In the first place we impose the
condition $\mathbf{F}(\mathbb{R}_{+}^{4})\subset\mathbb{R}_{+}^{4}$. This is
equivalent to making functions $f$ and $g$ verify $f([0,N^{1}]\times
\lbrack0,N^{2}])\subset\lbrack0,N^{1}]$ and $g([0,N^{1}]\times\lbrack
0,N^{2}])\subset\lbrack0,N^{2}]$, respectively, for any nonnegative $N^{1}$
and $N^{2}$:

\begin{itemize}
\item A first necessary condition is that $f(N^{1},N_{I}^{2})\geq0$, i.e.,
$(1-\gamma)N^{1}\geq0$, which implies that $\gamma\in(0,1]$.

\item We also need that $f(0,N^{2})\leq N^{1}$, i.e., $\dfrac{\beta N^{1}%
N^{2}}{N^{1}+N^{2}}\leq N^{1}$ for all $N^{1},N^{2}\geq0$. Since the limit of
the second term when $N^{2}$ tends to infinity is $\beta N^{1},$ we also
assume that $\beta\in(0,1]$.
\end{itemize}

These two necessary conditions turn out to be also sufficient for
$\mathbf{F}(\mathbb{R}_{+}^{4})\subset\mathbb{R}_{+}^{4}$. Indeed, let
\begin{equation}
\gamma,\beta\in(0,1] \label{condbg}%
\end{equation}
and $(N_{I}^{1},N_{I}^{2})\in\lbrack0,N^{1}]\times\lbrack0,N^{2}]$. Then
\[
f\left(  N_{I}^{1},N_{I}^{2}\right)  =(1-\gamma)N_{I}^{1}+\dfrac{\beta}%
{N^{1}+N^{2}}\left(  N^{1}-N_{I}^{1}\right)  \left(  N_{I}^{1}+N_{I}%
^{2}\right)  \geq0,
\]
and
\[
f\left(  N_{I}^{1},N_{I}^{2}\right)  \leq N_{I}^{1}+\beta\left(  N^{1}%
-N_{I}^{1}\right)  \leq N_{I}^{1}+\left(  N^{1}-N_{I}^{1}\right)  =N^{1},
\]
as we wanted. An analogous reasoning proves that also $g$ verifies the
required conditions.

Now, given the above conditions, it is immediate to check that $\mathbf{F}%
(\Omega)\subset\Omega$.

The asymptotic behaviour of the solutions of systems \eqref{fm} and \eqref{rm}
can be related by making use of results in \cite{Sanz08}
regarding approximate reduction techniques. The condition for the results to
hold is that $\mathbf{F}^{(k)}$ converges to a map $\bar{\mathbf{F}}$
uniformly on compact sets of $\Omega$ and the same happens with their
differentials, i.e., $\lim_{k\rightarrow\infty}D\mathbf{F}^{(k)}(N_{S}%
,N_{I})=D\bar{\mathbf{F}}(N_{S},N_{I})$ uniformly on compact sets.

To find map $\bar{\mathbf{F}}$ we study the long term behaviour of the
2-dimensional system describing the dynamics of infected individuals of both
species, $N_{I}^{1}$ and $N_{I}^{2}$:
\[
N_{I}^{1}(t+1)=f\left(  N_{I}^{1}(t),N_{I}^{2}(t)\right)  ,\ N_{I}%
^{2}(t+1)=g\left(  N_{I}^{1}(t),N_{I}^{2}(t)\right)
\]
Its equilibrium points are the disease free equilibrium $E_{0}=(0,0)$ and, if
$\beta>\gamma$, a positive equilibrium representing the disease endemicity
\[
E_{+}=\left(  (1-\frac{\gamma}{\beta})N^{1},(1-\frac{\gamma}{\beta}%
)N^{2}\right)
\]
that we could also express in terms of $\nu$ \eqref{nu}, or $R_{0}$,
\[
E_{+}=\left(  (1-\nu)N^{1},(1-\nu)N^{2}\right) =\left(  (1-\frac{1}{R_{0}%
})N^{1},(1-\frac{1}{R_{0}})N^{2}\right) .
\]
By linearization it is easy to prove that $E_{0}$ is A.S. if $\beta<\gamma$
and $E_{+}$ is A.S. if $\beta>\gamma$.

Henceforth we assume Hypothesis 1: $0<\gamma<\beta\leq1$.

Now, we define our candidate to be $\bar{\mathbf{F}}$ as
\[
\begin{array}{l}
\bar{\mathbf{F}}\left(  N_{S}^{1},N_{I}^{1},N_{S}^{2},N_{I}^{2}\right) =
 \\ \qquad \qquad  \left(  \nu(N_{S}^{1}+N_{I}^{1}),(1-\nu)(N_{S}^{1}+N_{I}^{1}),\nu(N_{S}
^{2}+N_{I}^{2}),(1-\nu)(N_{S}^{2}+N_{I}^{2})\right)  .
\end{array}
\]

We restrict out attention to the case in which there are infected individuals
in the initial population, i.e., we work in the set $\Omega$ \eqref{omega}
which is positively invariant for $\mathbf{F}$.

Next we prove the required convergence results of $\mathbf{F}^{(k)}$ to
$\bar{\mathbf{F}}$.

\begin{proposition}
\label{pcu} Let $\mathbf{F}$ be the map defined in \eqref{fd} with $\beta
_{ij}=\beta$ and $\gamma_{i}=\gamma$ for $i,j=1,2$, and let Hypothesis 1 hold.
Then the following two limits exist uniformly on compacts sets of $\Omega$ \eqref{omega}:

\begin{enumerate}
\item $\displaystyle \lim_{k\rightarrow\infty} \mathbf{F}^{(k)}(N^{1}_{S},
N^{1}_{I}, N^{2}_{S},N^{2}_{I})=\bar{\mathbf{F}}(N^{1}_{S}, N^{1}_{I},
N^{2}_{S},N^{2}_{I})$.

\item $\displaystyle \lim_{k\rightarrow\infty} D\mathbf{F}^{(k)}(N^{1}_{S},
N^{1}_{I}, N^{2}_{S},N^{2}_{I})=D\bar{\mathbf{F}}(N^{1}_{S}, N^{1}_{I},
N^{2}_{S},N^{2}_{I})$.
\end{enumerate}
\end{proposition}

\begin{proof}
In order to prove the result we define the sets
\begin{align*}
V_{i} &  :=\left\{  \left(  N_{S}^{1},N_{I}^{1},N_{S}^{2},N_{I}^{2}\right)
\in\Omega:N_{S}^{i}+N_{I}^{i}=0\right\}  ,\ i=1,2,\\
\hat{\Omega} &  :=\left\{  \left(  N_{S}^{1},N_{I}^{1},N_{S}^{2},N_{I}%
^{2}\right)  \in\Omega:N_{S}^{1}+N_{I}^{1}>0,N_{S}^{2}+N_{I}^{2}>0\right\}  .
\end{align*}
In order to carry out the proof we will show that $\lim_{k\rightarrow\infty
}\mathbf{F}^{(k)}=\bar{\mathbf{F}}$ and that\newline$\lim_{k\rightarrow\infty
}D\mathbf{F}^{(k)}=\bar{\mathbf{F}}$ uniformly on compact sets of $V_{1}$ (*),
on compact sets of $V_{2}$ (**) and on compact sets of $\hat{\Omega}$.

On $A_{1}$ (resp. $A_{2}$) map $\mathbf{F}$ has the form $\mathbf{F}%
(\mathbf{N})=\left(  0,0,\mathbf{\hat{F}}(N_{S}^{2},N_{I}^{2})\right)  $
(resp. $\mathbf{F}(\mathbf{N})=\left(  \mathbf{\hat{F}}(N_{S}^{1},N_{I}%
^{1}),0,0\right)  $) where
\[
\mathbf{\hat{F}}(x_{1},x_{2})=\left(  x_{1}-\beta\dfrac{x_{1}x_{2}}%
{x_{1}+x_{2}}+\gamma x_{1},x_{2}+\beta\dfrac{x_{1}x_{2}}{x_{1}+x_{2}}-\gamma
x_{1}\right)
\]
and therefore in order to prove (*) and (**) it suffices to show that
$\lim_{k\rightarrow\infty}\mathbf{\hat{F}}^{(k)}=\bar{\mathbf{F}}^{\ast}$ and
that $\lim_{k\rightarrow\infty}D\mathbf{\hat{F}}^{(k)}=\bar{\mathbf{F}}$ in
compact sets of $\left\{  \left(  x_{1},x_{2}\right)  :x_{2}>0,\ x_{1}%
+x_{2}>0\right\}  $. This result was proved in Lemma A.1 of
\cite{Bravo17b}.

Now we turn our attention to the uniform convergence in compact sets of
$\hat{\Omega}$. We begin by making a change of variables and expressing the
map $\mathbf{F}$ in terms of new variables: $x:=N_{I}^{1}+N_{I}^{2}$ the total
number of infected individuals, $y:=\dfrac{N_{I}^{1}}{N_{S}^{1}+N_{I}^{1}%
}-\dfrac{N_{I}^{2}}{N_{S}^{2}+N_{I}^{2}}$ the difference between the fraction
of infected individuals in the first and the second species, $z:=N^{1}%
=N_{S}^{1}+N_{I}^{1}$ the total population of the first species, and
$w:=N^{2}=N_{S}^{2}+N_{I}^{2}$ the total population of the second species.

Let $G$ denote the map associated to this change of variables:
\[
\begin{array}{l}
(x,y,z,w)=G(N_{S}^{1},N_{I}^{1},N_{S}^{2},N_{I}^{2}):=\\
\qquad \qquad  \qquad \qquad (N_{I}^{1}+N_{I}^{2},\dfrac{N_{I}^{1}}{N_{S}^{1}+N_{I}^{1}}
-\dfrac{N_{I}^{2}}{N_{S}^{2}+N_{I}^{2}},N_{S}^{1}+N_{I}^{1},N_{S}^{2}+N_{I}^{2}),
\end{array}
\]
and let $G^{-1}$ be its inverse map.

Map $\mathbf{F}$ in the new variables is easily found to be%
\[
\begin{array}{l}
\mathbf{H}(x,y,z,w):=G(\mathbf{F}(G^{-1}(x,y,z,w)))=\\
\qquad \qquad  \qquad  \qquad  \qquad  \left(  (1+\beta
-\gamma-\dfrac{\beta x}{z+w})x,(1-\gamma-\dfrac{\beta x}{z+w})y,z,w\right)  .
\end{array}
\]
Therefore we need to prove the uniform convergence on compact sets of
$G(\hat{\Omega})$ of its iterates $\mathbf{H}^{(k)}$ to map $\mathbf{\bar{F}}$
expressed in terms of the new variables, that we denote $\bar{\mathbf{H}}$:
\[
\bar{\mathbf{H}}(x,y,z,w):=G(\mathbf{\bar{F}}(G^{-1}(x,y,z,w)))=\left(
(1-\nu)(z+w),0,z,w\right)  .
\]
The uniform convergence of the differentials of $\mathbf{H}^{(k)}$ to the
differential of $\bar{\mathbf{H}}$ should also be proved.

With the help of the function
\[
\phi(x):=(1+\beta-\gamma-\beta x)x,
\]
we can express $\mathbf{H}$ in the following form
\[
\mathbf{H}(x,y,z,w)=\left(  (z+w)\phi(\frac{x}{z+w}),(1-\gamma-\frac{\beta
x}{z+w})y,z,w\right)  ,
\]
and its $k$-th iterate as
\[
\mathbf{H}^{(k)}(x,y,z,w)=\left(  (z+w)\phi^{(k)}(\frac{x}{z+w}),y\prod
_{i=0}^{k-1}\left(  1-\gamma-\beta\phi^{(i)}(\frac{x}{z+w})\right)
,z,w\right)  .
\]
Note that in $G(\hat{\Omega})$ we have $0<x\leq z+w$. \ Let $K\subset
G(\Omega)$ be a compact set. Then there exist numbers
\begin{equation}
\displaystyle M_{K}=\max_{K}(z+w),\ \displaystyle a_{K}=\min_{K}(\frac{x}%
{z+w})>0\text{ and }\displaystyle b_{K}=(\max_{K}\frac{x}{z+w})\leq1.
\label{a1}%
\end{equation}

Function $\phi$ verifies that $\phi((0,1])\subset(0,1]$ and so
\begin{equation}
\left\vert 1-\gamma-\beta\phi^{(i)}(\frac{x}{z+w})\right\vert \leq
\max\{1-\gamma,|1-\gamma-\beta|\}=:c<1, \label{a2}%
\end{equation}
where in the last inequality we have used Hypothesis 1. Then we have
\[%
\begin{array}
[c]{c}%
\displaystyle\max_{K}\Vert\mathbf{H}^{k}(x,y,z,w)-\bar{\mathbf{H}%
}(x,y,z,w)\Vert_{1}=\\
\displaystyle\max_{K}\left(  |(z+w)(\phi^{(k)}(\frac{x}{z+w})-(1-\nu
))|+|y\prod_{i=0}^{k-1}\left(  1-\gamma-\beta\phi^{(i)}(\frac{x}{z+w})\right)
|\right)  \leq\\
\displaystyle M_{K}\max_{K}\left(  \phi^{(k)}(\frac{x}{z+w})-(1-\nu)\right)
+c^{k}.
\end{array}
\]
where we have used that one has $\left\vert y\right\vert \leq1$ in
$G(\hat{\Omega})$. Now it is straightforward to see that the solutions of the
scalar difference equation $x(t+1)=\phi(x(t))$ with initial conditions
$x(0)\in\lbrack a_{K},b_{K}]\subset(0,1]$ converge monotonically and,
therefore, uniformly on $[a_{K},b_{K}]$, to $1-\nu$. Therefore we have proved
the uniform convergence of $\mathbf{H}^{(k)}$ to $\bar{\mathbf{H}}$ on compact
sets of $G(\hat{\Omega})$.

To prove the uniform convergence of the differential of $\mathbf{H}^{(k)}$ we
start by expressing $D\mathbf{H}^{(k)}$ in terms of the derivatives of
$\phi^{k}$:
\begin{equation}
D\mathbf{H}^{(k)}(x,y,z,w)=\left(
\begin{array}
[c]{cccc}%
(\phi^{(k)})^{\prime}(\dfrac{x}{z+w}) & 0 & d_{13} & d_{14}\\
\rule{0ex}{3ex}d_{21} & \prod_{i=0}^{k-1}\left(  1-\gamma-\beta\phi
^{(i)}(\dfrac{x}{z+w})\right)  & d_{23} & d_{24}\\
\rule{0ex}{3ex}0 & 0 & 1 & 0\\
0 & 0 & 0 & 1
\end{array}
\right)  , \label{zzz}%
\end{equation}
with
\[
d_{13}=d_{14}=\phi^{(k)}(\frac{x}{z+w})-\frac{x}{z+w}(\phi^{(k)})^{\prime
}(\frac{x}{z+w}),
\]%
\[
d_{21}=y\sum_{j=0}^{k-1}\left(  -\frac{\beta}{z+w}(\phi^{(j)})^{\prime}%
(\dfrac{x}{z+w})\prod_{{\scriptsize i=0,\ i\neq j}}^{k-1}\left(
1-\gamma-\beta\phi^{(i)}(\frac{x}{z+w})\right)  \right)  ,
\]
and
\[
d_{23}=d_{24}=y\sum_{j=0}^{k-1}\left(  \frac{\beta x}{(z+w)^{2}}(\phi
^{(j)})^{\prime}(\dfrac{x}{z+w})\prod_{{\scriptsize i=0,\ i\neq j}}%
^{k-1}\left(  1-\gamma-\beta\phi^{(i)}(\frac{x}{z+w})\right)  \right)  .
\]
Now we want to show that the following limit is uniform on compact sets of
$(0,1]$
\[
\lim_{k\rightarrow\infty}(\phi^{k})^{\prime}(x)=0.
\]
Since $|\phi^{\prime}(1-\nu)|<1,$ there exist $\alpha<1$ and a neighbourhood
$I\subset(0,1]$ of $1-\nu$ such that for every $x\in I$ we have $|\phi
^{\prime}(x)|<\alpha$. The uniform convergence to $1-\nu$ of the solutions of
the scalar difference equation $x(t+1)=\phi(x(t))$ with initial conditions
$x_{0}\in(0,1]$ together with the chain rule to obtain $(\phi^{k})^{\prime
}(x)=\prod_{{\scriptsize i=0}}^{k-1}\phi^{\prime}\left(  \phi^{(i)}(x)\right)
$ yield the result.

This last result, together with \eqref{a1}, \eqref{a2} and \eqref{zzz}
straightforwardly imply the uniform convergence on compact sets of
$D\mathbf{H}^{(k)}(x,y,z,w)$ to
\[
D\bar{\mathbf{H}}(x,y,z,w)=\left(
\begin{array}
[c]{cccc}%
0 & 0 & 1-\nu & 1-\nu\\
0 & 0 & 0 & 0\\
0 & 0 & 1 & 0\\
0 & 0 & 0 & 1
\end{array}
\right)  .
\]

\end{proof}

The next theorem relates the asymptotic behavior of systems (\ref{fm}) and
(\ref{rm}) for large enough values of parameter $k$ when the hypotheses of
Proposition \ref{pcu} hold.

\begin{theorem}
\label{tha} Let the hypotheses of Proposition \ref{pcu} hold. Let $(N^{1\ast
},N^{2\ast})$ be a hyperbolic equilibrium point of (\ref{rm}). Then there
exists $k_{0}\in\mathbb{N}$ such that for each $k\geq k_{0}$ there exists a
hyperbolic equilibrium point $\left(  N_{S,k}^{1\ast},N_{I,k}^{1\ast}%
,N_{S,k}^{2\ast},N_{I,k}^{2\ast}\right)  $ of (\ref{fm}) satisfying
\[
\lim_{k\rightarrow\infty}\left(  N_{S,k}^{1\ast},N_{I,k}^{1\ast}%
,N_{S,k}^{2\ast},N_{I,k}^{2\ast}\right)  =\left(  \nu N^{1\ast},(1-\nu
)N^{1\ast},\nu N^{2\ast},(1-\nu)N^{2\ast}\right)  .
\]
Moreover, let $k\geq k_{0}$ be fixed:

\begin{enumerate}
\item If $(N^{1\ast},N^{2\ast})$ is asymptotically stable (resp. unstable)
then \\$\left(  N_{S,k}^{1\ast},N_{I,k}^{1\ast},N_{S,k}^{2\ast},N_{I,k}^{2\ast
}\right)  $ is asymptotically stable (resp. unstable).
\item In the case of $(N^{1\ast},N^{2\ast})$ being asymptotically stable, \\
if $\left(  N_{S}^{1}(0)+N_{I}^{1}(0),N_{S}^{2}(0)+N_{I}^{2}(0)\right)$ is in the basin of attraction of $(N^{1\ast},N^{2\ast}),$ then $\left(  N_{S}%
^{1}(0),N_{I}^{1}(0),N_{S}^{2}(0),N_{I}^{2}(0)\right)$ is in the basin of
attraction of \\
$\left(  N_{S,k}^{1\ast},N_{I,k}^{1\ast},N_{S,k}^{2\ast},N_{I,k}^{2\ast}\right)$.
\end{enumerate}
Analogous results hold for periodic solutions.
\end{theorem}

\begin{proof}
It is a direct consequence of the results in \cite{Sanz08} and
Proposition \ref{pcu}.
\end{proof}

\section{Dynamics of the reduced system}

\label{app2}

The proof of Lemma \ref{prop:lema1}, Proposition \ref{prop:prop01} and
Theorems \ref{prop:prop2} and \ref{prop:prop3} follows reasonings very similar
to those of Lemma 2, Propositions 3 and 5, and Theorems 4 and 6 in
\cite{Bravo17a}. In that reference, the system under analysis (Eq.
(20) in that reference) is a particular case of system (\ref{sist})
corresponding to making $r_{I}^{2}=0$. As a consequence, here both isoclines
are hyperbolas whereas in \cite{Bravo17a} $\Gamma_{1}$ is a hyperbola
and $\Gamma_{2}$ is a straight line. Therefore, in the proofs below we
concentrate on the differences with \cite{Bravo17a} and refer the
reader to that reference for those reasonings that are identical.

\noindent\textbf{Proof of Lemma }\ref{prop:lema1}. This is essentially Lemma 2
in \cite{Bravo17a} except for the fact that here we apply the
reasonings therein to both $S_{1}$ and $S_{2}$ and not only to $S_{1}$. We
also need to give conditions for the isoclines $S_{i}$ to be degenerate which
are not studied in that reference. From the classical theory of conics one has
that the two discriminants for $S_{i}$ are%
\small
\begin{align*}
\delta_{i}  &  :=\det\left(
\begin{array}
[c]{cc}%
c_{S1}^{i}c_{I1}^{i} & \frac{1}{2}\left(  c_{S1}^{i}c_{I2}^{i}+c_{S2}%
^{i}c_{I1}^{i}\right) \\
\frac{1}{2}\left(  c_{S1}^{i}c_{I2}^{i}+c_{S2}^{i}c_{I1}^{i}\right)  &
c_{S2}^{i}c_{I2}^{i}%
\end{array}
\right)  =-\frac{1}{4}\left(  c_{S1}^{i}c_{I2}^{i}-c_{S2}^{i}c_{I1}%
^{i}\right)  ^{2}\leq0\\
\Delta_{i}  &  :=\det\left(
\begin{array}
[c]{ccc}%
c_{S1}^{i}c_{I1}^{i} & \frac{c_{S1}^{i}c_{I2}^{i}+c_{S2}^{i}c_{I1}^{i}}{2} &
-\frac{c_{S1}^{i}(r_{I}^{i}-1)+c_{I1}^{i}(r_{S}^{i}-1)}{2}\\
\frac{c_{S1}^{i}c_{I2}^{i}+c_{S2}^{i}c_{I1}^{i}}{2} & c_{S2}^{i}c_{I2}^{i} &
-\frac{c_{S2}^{i}(r_{I}^{i}-1)+c_{I2}^{i}(r_{S}^{i}-1)}{2}\\
-\frac{c_{S1}^{i}(r_{I}^{i}-1)+c_{I1}^{i}(r_{S}^{i}-1)}{2} & -\frac{c_{S2}%
^{i}(r_{I}^{i}-1)+c_{I2}^{i}(r_{S}^{i}-1)}{2} & 1-r_{S}^{i}-r_{I}^{i}%
\end{array}
\right) \\
&  =\frac{1}{4}r_{S}^{i}r_{I}^{i}\left(  c_{S1}^{i}c_{I2}^{i}-c_{S2}^{i}%
c_{I1}^{i}\right)  ^{2}\geq0
\end{align*}%
\normalsize
It is well known that hyperbola $S_{i}$ is degenerate if and only if
$\delta_{i}=0$, i.e., if and only if $c_{S1}^{i}c_{I2}^{i}=c_{S2}^{i}%
c_{I1}^{i}$ and, in that case, $\Delta_{i}=0$ and therefore it corresponds to
two parallel lines.

\noindent\textbf{Proof of Proposition }\ref{prop:prop01}. The proof of (a) is
straightforward. Regarding (b), it is immediate to check that, for
$i,j\in\{1,2\}$, $i\neq j$, the map $H_{i}(x_{1},x_{2})$ is strictly
increasing as a function of $x_{i}$ and strictly decreasing as a function of
$x_{j}$. Therefore it follows trivially that if $x,\ x^{\prime}\in
\mathbf{R}_{+}^{2}$ are distinct points with $x\leq_{K}x^{\prime}$ then
$H(x)<_{K}H(x^{\prime})$.

(c) Let us consider map $H:\bar{S}\rightarrow\bar{S},$ where $S$ is defined in
(\ref{eqS}):

(i) It is immediate to check that $\bar{S}$ contains order intervals \cite[p.
345]{Smith98} and is $\leq_{K}$-convex \cite[p. 339]{Smith98}.

(ii) In order to show that $\det DH(x_{1},x_{2})>0$ for $(x_{1},x_{2})\in
\bar{S}$, where $D$ denotes differential, we have made use of Matlab Symbolic
Math Toolbox, and obtained that the resulting expression is the product of a
number of factors all of which are strictly positive in $\bar{S}$.

(iii) Direct calculations prove that $DH(x_{1},x_{2})$ is $K$-positive in
$\bar{S}$ \cite[p. 338]{Smith98}, i.e., $\dfrac{\partial H_{1}%
}{\partial x_{1}}>0$, $\dfrac{\partial H_{2}}{\partial x_{2}}>0$,
$\dfrac{\partial H_{1}}{\partial x_{2}}\leq0$ and $\dfrac{\partial H_{1}%
}{\partial x_{1}}\leq0$ en $\bar{S}$.

(iv) Finally, it is immediate to check that $\bar{S}$ is compact and connected
and that $H^{-1}(0,0)$ is a single point.

Using the properties (i) through (iv) above, \cite[Lemma 4.3]{Smith98}
guarantees that $H$ verifies property (O+) \cite[p. 343]{Smith98}.
Now, using (a) and (b) and property (O+) we can apply Theorem 4.2 in
\cite{Smith98} so that all orbits in $\mathbf{R}_{+}^{2}$ are
eventually componentwise monotone and converge to an equilibrium.

\noindent\textbf{Proof of Theorem }\ref{prop:prop2}. It corresponds to
Proposition 5 in \cite{Bravo17a}. The reasonings to carry out its
proof are literally identical to the ones used in that reference.

\noindent\textbf{Proof of Theorem }\ref{prop:prop3}. Part 1 of this result
corresponds to Proposition 5 in \cite{Bravo17a} and the reasonings
therein can be translated literally to this case. Part 2 is immediate taking
into account that the positive equilibria are the intersections of $\Gamma
_{1}$ and $\Gamma_{2}$ and both are the graphs of strictly decreasing
functions of $x_{1}$.

Let us consider part 3, in which the eight different cases described in
(\ref{casos}) are studied. To start with, we point out that Lemmas A.1. and
A.2. in \cite{Bravo17a} translate directly to our setting. Indeed,
the first one hinges in the properties proved in Proposition \ref{prop:prop01}%
, whereas the second is a consequence of the facts that $\dfrac{\partial
H_{1}}{\partial x_{2}}<0$ and $\dfrac{\partial H_{1}}{\partial x_{1}}<0$ in
$\mathring{\mathbf{R}}_{+}^{2}$ and that $\det DH(x_{1},x_{2})>0$ for
$(x_{1},x_{2})\in\bar{S}$. Using Lemma A.2. and the fact that isoclines are
not tangent at equilibria, it follows that all the positive equilibria are hyperbolic.

Now, the statement regarding our scenarios \textbf{A}$_{1}$, \textbf{B}$_{1}$,
\textbf{C}$_{0}$, \textbf{C}$_{2}$ and \textbf{D}$_{0}$, correspond,
respectively, to the statements for cases \textbf{A}, \textbf{B},
\textbf{C}$_{1a}$, \textbf{C}$_{1b}$ and \textbf{C}$_{2}$ in Theorem 6\ in
\cite{Bravo17a}, the only difference being that in that reference
isocline $\Gamma_{2}$ is a straight line instead of a hyperbola. This does not
alter in any way the reasonings carried out therein and so they translate
literally to our setting.

Regarding our case \textbf{D}$_{2}$, it is the reciprocal of case
\textbf{C}$_{2}$ interchanging $x_{1}$ with $x_{2}$, $E_{1}^{\ast}$ with
$E_{2}^{\ast}$, and $E_{3}^{\ast}$ with $E_{4}^{\ast}$, and so the proof for
this case follows from the one corresponding to \textbf{C}$_{2}$.

Let us turn prove the statements regarding the only remaining cases,
\textbf{A}$_{3}$ and \textbf{B}$_{3}$, in which there are three positive
equilibria. Let us consider the six open connected regions in which
$\Gamma_{1}$ and $\Gamma_{2}$ divide $\mathbb{R}_{+}^{2}$ in cases (see Figure
\ref{fig1}) \textbf{A}$_{3}$ and \textbf{B}$_{3}$. Let us define $U_{0}$ as
the region whose adherence contains the origin, $U_{\infty}$ as the only
unbounded region, $U_{23}$ as the region limited by $\Gamma_{1}$ and
$\Gamma_{2}$ whose adherence contains equilibria $E_{2}^{\ast}$ and
$E_{3j}^{\ast}$ and similarly for $U_{34}$, $U_{45}$ and $U_{51}$. We have
then that $\mathbb{R}_{+}^{2}=U_{0}\cup U_{\infty}\cup\bar{U}_{23}\cup\bar
{U}_{34}\cup\bar{U}_{45}\cup\bar{U}_{51}$. We can divide $\mathbb{R}_{+}^{2}$
into the four disjoint sets%
\begin{align*}
W  &  :=\left\{  x\in\mathbf{R}_{+}^{2}:H(x)\leq_{K}x\right\}  ,\ W^{\prime
}:=\left\{  x\in\mathbf{R}_{+}^{2}:H(x)\geq_{K}x\right\} \\
T  &  :=\left\{  x\in\mathbf{R}_{+}^{2}:H(x)<x\right\}  ,\ T^{\prime
}:=\left\{  x\in\mathbf{R}_{+}^{2}:H(x)>x\right\}  .
\end{align*}

Let us consider case \textbf{A}$_{3}$. Using the monotonicity of the two
components of map $H$ (see Figure \ref{fig1}), we have that in this case
$W=\bar{U}_{34}\cup\bar{U}_{51}$, $W^{\prime}=\bar{U}_{23}\cup\bar{U}_{45}$,
$T=U_{\infty}$ and $T^{\prime}=U_{0}$. Using Lemma A.1. in
\cite{Bravo17a} it follows that $\bar{U}_{34}\cup\bar{U}_{51}$ and
$\bar{U}_{23}\cup\bar{U}_{45}$ are forward invariant for $H$, that orbits
starting in $U_{0}$ can not enter $U_{\infty}$ and that orbits starting in
$U_{\infty}$ can not enter $U_{0}$. Since $\bar{U}_{23}\cup\bar{U}_{45}$ is
forward invariant, $E_{3}^{\ast}\in\bar{U}_{23}$ is a fixed point and $\bar
{U}_{23}$ is connected, $H(\bar{U}_{23})$ can not intersect $\bar{U}_{45}$ and
therefore $\bar{U}_{23}$ must be forward invariant. A similar reasoning proves
that $\bar{U}_{45}$, $\bar{U}_{34}$ and $\bar{U}_{51}$ are forward invariant.

By the monotonicity of $H$ in the different regions we have that orbits
starting in $U_{34}$, $U_{45}$, $U_{23}$ and $U_{51}$ respectively, can not
converge to $E_{4}^{\ast}$, $E_{4}^{\ast}$, $E_{2}^{\ast}$ and $E_{1}^{\ast}$,
respectively and using that each orbit must converge to a fixed point we have
that they must converge, respectively, to $E_{3}^{\ast}$, $E_{5}^{\ast}$,
$E_{3}^{\ast}$ and $E_{5}^{\ast}$. In particular $E_{4}^{\ast}$ is unstable.
Let us now show that $E_{3}^{\ast}$ is attracting. Let $\varepsilon>0$ be
small enough, let $B(E_{3}^{\ast},\varepsilon)$ be the open ball with center
$E_{3}^{\ast}$ and radius $\varepsilon$ and let $x\in B(E_{3}^{\ast
},\varepsilon)$, so that $x$ must belong to one and only one of the sets
$\bar{U}_{23}$, $\bar{U}_{34}$, $U_{0}$ and $U_{\infty}$. Let us consider in
turn the three following possibilities: (i) $x\in\bar{U}_{23}\cup\bar{U}_{34}%
$, (ii) $x\in U_{0}$ and (iii) $x\in U_{\infty}$.

In scenario (i) we have already shown that the corresponding orbit converges
to $E_{3}^{\ast}$. (ii) We have proved that the corresponding orbit
$H^{(n)}(x)$ can not enter $U_{\infty}$ and so either $H^{(n)}(x)$ enters
$\bar{U}_{23}\cup\bar{U}_{34}$ for a certain $n,$ and in that case we already
know that the orbit must converge to $E_{3}^{\ast}$ or, on the contrary,
$H^{(n)}(x)$ remains in $U_{0}$ for all $n\geq0$. In this latter case the fact
that $x$ is $\varepsilon$-close to $E_{3}^{\ast}$ and the monotonicity in this
region precludes convergence of the orbit to $E_{1}^{\ast}$, $E_{2}^{\ast}$,
$E_{3}^{\ast}$ or to $E_{5}^{\ast},$ and so it must necessarily converge to
$E_{3}^{\ast}$. (iii) In this case we can carry out a reasoning completely
analogous to that of (ii). As a conclusion we have that $E_{3}^{\ast}$
attracts the open ball $B(E_{3}^{\ast},\varepsilon)$ and so it is attracting.
The proof that $E_{5}^{\ast}$ is attracting is carried out similarly. In order
to show that $E_{4}^{\ast}$ is a saddle cause it attracts points different
from itself, and that it can not attract any open set, we can use the same
reasonings carried out for $E_{4}^{\ast}$ in case \textbf{C}$_{1b}$ of Theorem
6 in \cite{Bravo17a}. Finally, the monotonicity of $H$ in the
neighbourhood of $E_{1}^{\ast}$ and of $E_{2}^{\ast}$ precludes the
convergence of any orbit to these two points and, therefore, except for the
orbits that converge to the saddle $E_{4}^{\ast}$, the rest of the orbits in
$\mathring{\mathbf{R}}_{+}^{2}$ converge to $E_{3}^{\ast}$ or to $E_{5}^{\ast
}$.

Regarding the statements for case \textbf{B}$_{3}$, they are proved
analogously to those of case \textbf{A}$_{3}$.

\end{document}